\documentclass[11pt]{article}

\usepackage{fullwidth}
\usepackage{color}
\usepackage[color,matrix,arrow]{xy}
\usepackage[top=1.5in, bottom=1.5in, left=1in, right=1in]{geometry}
\usepackage{amsgen}
\usepackage{amsmath}
\usepackage{amstext}
\usepackage{amsbsy}
\usepackage{amsopn}
\usepackage{amsfonts}
\usepackage{amssymb}
\usepackage{eepic}
\usepackage{graphicx}
\usepackage{epsf}
\usepackage{pstricks}
\xyoption{all}

\def\Box{\square}

\def\tra#1{\smash{\mathop{\mid\kern
-1pt\joinrel\relbar\joinrel\relbar}\limits^{*}_{#1}}}
\def\longtra#1{\smash{\mathop{\mid\kern
-1pt\joinrel\relbar\joinrel\relbar\joinrel\relbar}\limits^{*}_{#1}}}
\def\vlongtra#1{\smash{\mathop{\mid\kern
-1pt\joinrel\relbar\joinrel\relbar\joinrel\relbar\joinrel\relbar}\limits^{*}_{#1}}}
\def\vvlongtra#1{\smash{\mathop{\mid\kern
-1pt\joinrel\relbar\joinrel\relbar\joinrel\relbar\joinrel\relbar\joinrel\relbar}\limits^{*}_{#1}}}
\def\vvvlongtra#1{\smash{\mathop{\mid\kern
-1pt\joinrel\relbar\joinrel\relbar\joinrel\relbar\joinrel\relbar\joinrel\relbar\joinrel\relbar}\limits^{*}_{#1}}}
\def\etra#1{\smash{\mathop{\mid\kern
-1pt\joinrel\relbar\joinrel\relbar}\limits_{#1}}}

\def\iff{\Leftrightarrow}

\def\B{{\cal{B}}}

\def\F{{\cal{F}}}

\def\L{{\cal{L}}}

\def\P{{\cal{P}}}

\def\H{{\mathcal{H}}}

\def\bi{\begin{itemize}}
\def\ei{\end{itemize}}
\def\beq{\begin{equation}}
\def\eeq{\end{equation}}

\newtheorem{T}{Theorem}[section]
\newcommand{\bt}{\begin{T}}
\newcommand{\et}{\end{T}}
\newcommand{\ftd}{$\square$\end{T}}

\newtheorem{Proposition}[T]{Proposition}
\newcommand{\bp}{\begin{Proposition}}
\newcommand{\ep}{\end{Proposition}}
\newcommand{\fpd}{$\square$\end{Proposition}}

\newtheorem{Lemma}[T]{Lemma}
\newcommand{\bl}{\begin{Lemma}}
\newcommand{\el}{\end{Lemma}}
\newcommand{\fld}{$\square$\end{Lemma}}

\newtheorem{Corol}[T]{Corollary}
\newcommand{\bc}{\begin{Corol}}
\newcommand{\ec}{\end{Corol}}
\newcommand{\fcd}{$\square$\end{Corol}}

\newtheorem{Result}[T]{Result}
\newcommand{\br}{\begin{Result}}
\newcommand{\er}{\end{Result}}
\newcommand{\frd}{$\square$\end{Result}}

\newtheorem{Example}[T]{Example}
\newcommand{\be}{\begin{Example}}
\newcommand{\ee}{\end{Example}}

\newtheorem{Problem}[T]{Problem}
\newcommand{\bq}{\begin{Problem}}
\newcommand{\eq}{\end{Problem}}

\newtheorem{Remark}[T]{Remark}
\newcommand{\brm}{\begin{Remark}}
\newcommand{\erm}{\end{Remark}}

\newtheorem{Definition}[T]{Definition}
\newcommand{\bd}{\begin{Definition}}
\newcommand{\ed}{\end{Definition}}

\newcommand{\proof}
   {\par\medbreak\noindent{\bf Proof}.\enspace}

\newcommand{\qed}{
$\Box$
\par\bigbreak}

\normalbaselines
\def\abstract#1{\par\bigskip
\begingroup\small
\baselineskip=12truept
\begin{center}ABSTRACT\end{center}
\par\medskip\par\noindent
\null\hfill\hbox{\vbox{\hsize=5truein\noindent#1}}
\hfill\null\par\endgroup\par}

\title{On The Wilson Monoid of a Pairwise Balanced Design}
\author{{\bf Stuart Margolis, John Rhodes and Pedro V. Silva}}

\date{\today}

\begin{document}
\maketitle

\begin{center}\small

\bigskip

\end{center}

\abstract{We give a new perspective of the relationship between simple matroids of rank 3 and pairwise balanced designs, connecting Wilson's theorems and tools with the theory of truncated boolean representable simplicial complexes. We also introduce the concept of Wilson monoid $W(X)$ of a pairwise balanced design $X$. We present some general algebraic properties and study in detail the cases of Steiner triple systems up to 19 points, as well as the case where a single block has more than 2 elements.
}

\section{Introduction}

The purpose of this paper is to describe the deep connection between the theory of truncated boolean representable simplicial complexes ({\em TBRSC}) \cite{brsc, tbrsc} and R. Wilson's famous theorem that pairwise block designs ({\em PBDs}) exist for large enough sets meeting the usual necessary conditions on their parameters \cite{Wilson1, Wilson2, Wilson3}. In addition, we begin the algebraic study of the monoid of Wilson morphisms from a $PBD$ to itself. This gives important connections between the theory of matroids, truncated boolean representable simplicial complexes, design theory and semigroup theory that are mutually beneficial for all these fields. We outline these connections briefly in this section.

In matroid theory, the inverse operation of truncation is called erection \cite{Oxley}. The study of this operator 
was initiated by Crapo \cite{Crapo1} and plays an important part in matroid theory \cite{Knuth1, Nguyen1, Pendavingh1}.  In this paper we study erections for matroids of rank 3 within the context of the theory of $TBRSC$. That is, if $M$ is a matroid of rank 3, on the set of points $V$ and independent sets $\mathcal{H}$, then we wish
to compute the maximal boolean representable simplicial complex ({\em BRSC}) \cite{brsc} $M^{\varepsilon}$ whose truncation to rank 3 is $M$. Remarkably, this question is directly related to the fundamental papers of Wilson \cite{Wilson1,Wilson2,Wilson3}, one of the most important works in Combinatorics in the last 40 years. Indeed, we will see that the subsystems of a Pairwise Balanced Design ({\em PBDs}), (called flats of a $PBD$ by Wilson in \cite{Wilson1}) are precisely the flats of $M^{\varepsilon}$ in the sense of the theory of boolean representable simplicial complexes ({\em BRSC}) \cite{brsc}. As the lattice of subsystems of a {\em PBD} is rarely a geometric lattice, the connection to the work of Wilson was not studied by matroid theorists. It is only through the theory of {\em BRSC} that the lattice of subsystems of a $PBD$ plays its proper role.

Crapo \cite{Crapo1} proved that for matroids of any rank, the collection of all matroid erections form a lattice for the weak order and described the largest element in the lattice called the free erection. See also, \cite{Knuth1,Nguyen1,Pendavingh1}. 
The relationship between free erections and the {\em BRSC} $M^{\varepsilon}$ is discussed in \cite{tbrsc}.

More precisely, we recall that a {\em PBD} with parameter $\lambda = 1$ is well known to be equivalent, except for some trivial cases, to a matroid of rank 3 \cite{craporota}. Every such matroid $M$ defines a largest {\em BRSC}, $M^{\varepsilon}$, such that the truncation of $M^{\varepsilon}$ to rank 3 is $M$. An important result of this paper is to show that $M^{\varepsilon}$ is the {\em BRSC} defined by taking the {\em subsystems} of the corresponding {\em PBD} as the defining flats in the sense of \cite{brsc}. Recall that a subsystem of a {\em PBD} is a subset $X$ of the base set such that for each pair of distinct points $x,y \in X$, the unique block of the {\em PBD} containing them is contained in $X$. This latter {\em BRSC} is not in general a matroid itself, but gives us a way to lift the original matroid $M$ to higher dimensions. We give a number of illuminating examples.

In \cite{Wilson1}, Wilson defines a notion of morphism between $PBDs$ that gives the collection of $PBDs$ the structure of a category. In particular, if $X$ is a $PBD$, then the collection $W(X)$ of all morphisms from $X$ to itself forms a monoid that we call the Wilson monoid of $X$. Morphisms play a central role in \cite{Wilson1} where they are used to prove \cite[Section 8]{Wilson1}  what is now called Wilson's Fundamental Construction \cite[IV.2.1 Theorem 2.5]{handdesign}. This is one of the most important recursive constructions in design theory. Wilson's proof is by understanding the structure of the kernel of a morphism \cite[Theorem 8.1-8.2]{Wilson1}.

Despite their importance in Wilson's seminal work, morphisms were not subsequently developed. In particular, there has been no study of the algebraic properties of the monoid $W(X)$ of a $PBD$ $X$ and its relationship to the combinatorial and geometric properties of $X$. A major portion of this paper is devoted to developing these connections. We show that $W(X)$ consists precisely of the continuous partial functions on the collection of open sets $\mathcal{O}(X)$, in that inverse images of open sets are open. An open set is the complement of a subsystem of $X$. $\mathcal{O}(X)$ is closed under unions, but not necessarily intersection and thus we are working with a generalized version of a topology. Algebraically, $W(X)$ has a unique 0-minimal ideal $I(X)$ on which $W(X)$ acts faithfully on both the left and the right and is the largest monoid with this property. The connection between incidence structures and maximal faithful ideal extensions was studied by Dinitz and Margolis \cite{bibdtrans,bibdcont,projcont} in the 1980s.

As illuminating examples we look at two ``minimal" cases. First we look at Wilson monoids of Steiner triple systems, that is, $PBDs$ all of whose block sizes are 3. Then we look at $PBDs$ that have at most one block of size greater than 2. The case of triple systems indicates that almost every Wilson monoid consists of a unique 0-minimal ideal and its group of units, which is the automorphism group of the design. Monoids with this property are called small monoids. On the other hand, the triple system associated to both  affine $n$-space and projective $n$-space over the field of order 2 have monoids that contain the monoid of all $n \times n$ matrices over this field. This dichotomy between $PBDs$ that are ``small" and those that are ``big" is an important part of the theory. Thus, most triple systems are ``weeds'', in that they have no non-trivial automorphisms nor subsystems. Almost all of these have small monoids as Wilson monoid. On the other hand, the ``jewels'' are both rare and have a very intricate Wilson monoid.

%

%
%
%

\section{Boolean Representable Simplicial Complexes}

We review the basics of the theory of boolean representable simplicial complexes in this section. The reader is referred to \cite{brsc}. All the results mentioned here will be used throughout the paper without further reference.

All lattices and simplicial complexes in this paper are assumed to be
finite. Given a set $V$ and $n \geq 0$, we denote by $P_n(V)$
(respectively $P_{\leq n}(V)$)
the set
of all subsets of $V$ with precisely (respectively at most)
$n$ elements. 

A (finite) simplicial complex is a structure of the form $S =
(V,\H)$, where $V$ is a finite nonempty set and $\H \subseteq 2^V$
is nonempty and closed under taking subsets. 
The elements of $V$ and $\H$ are called respectively {\em points} and {\em independent sets}. 

A maximal independent set is called a {\em basis}. The maximum size of a basis is the {\em rank} of $S$. We say that $S$ is {\em pure} if all its bases have the same size. We say that $S =
(V,\H)$ is {\em simple} if $P_2(V) \subseteq \H$. 

A simplicial complex $M = (V,\H)$ is called a {\em matroid} if it
satisfies the {\em exchange property}:
\bi
\item[(EP)]
For all $I,J \in \H$ with $|I| = |J|+1$, there exists some
  $i \in I\setminus J$ such that $J \cup \{ i \} \in \H$.
\ei

There are many cryptomorphic definitions of matroids \cite{Oxley}. In this paper, since we are concerned with simplicial complexes with various notions of independence, we will always refer to a matroid via its simplicial complex of independent sets as above.

An important example of matroids are the {\em uniform matroids} $U_{k,n}$: for all $1 \leq k \leq n$, we write $U_{k,n} = (V,P_{\leq k}(V))$ with $|V| = n$.

A subset $\F$ of $2^V$ is called a {\em Moore family} if $V \in \F$ and $\F$ is closed under intersection (that is, a Moore family is a submonoid of the monoid of all subsets of $V$ under intersection). Every Moore family, under inclusion, constitutes a lattice (with intersection as meet and the determined join). We say that $X \subseteq V$ is a {\em transversal of the
successive differences} for a chain
$$F_0 \subset F_1 \subset \ldots \subset F_k$$
in $\F$ if $X$ admits an enumeration $x_1,\ldots , x_k$ such that $x_i \in F_i
\setminus F_{i-1}$ for $i = 1,\ldots,k$. 
We denote by ${\rm Tr}(\F)$ the set of transversals of the
successive differences for chains in $\F$.

We say that a simplicial complex $S = (V,\H)$ is {\em boolean representable} (BRSC) if $\H\; = {\rm Tr}(\F)$ for some Moore family $\F\, \subseteq 2^V$. Moreover, every BRSC can be obtained this way by taking as Moore family its {\em lattice of flats} (see \cite[Chapters 5 and 6]{brsc}):

We say that $X
\subseteq V$ is a {\em flat} of $S$ if
$$\forall I \in \H \cap 2^X \hspace{.2cm} \forall p \in V \setminus X
\hspace{.5cm} I \cup \{ p \} \in \H.$$
The set of all flats of $S$ is denoted by 
$L(S)$. Note that $V, \emptyset \in L(S)$ in all cases, and $L(S)$ is indeed a Moore family.

It follows from \cite[Corollary 5.2.7]{brsc} that a simplicial complex $S = (V,\H)$ is boolean
representable If and only if $\H\, = {\rm Tr}(L(S))$. Furthermore, the lattice $L(S)$ induces a closure operator on $2^V$ defined by
$$Cl(X) = \cap\{ F \in L(S) \mid X \subseteq F \}$$
for every $X \subseteq V$. 

An alternative characterization of BRSC is provided by boolean matrices \cite{brsc}, which explains the terminology.

All matroids are boolean representable \cite[Theorem
  5.2.10]{brsc}, but the converse is not true. Indeed, all matroids are pure but BRSC need not to be so. Unlike simple matroids, simple BRSC do not need to have a geometric lattice of flats \cite[Example 5.2.11]{brsc}.



\section{Truncated Boolean Representable Simplicial Complexes}
\label{stru}

Given a simplicial complex $S = (V,\mathcal{H})$ and $k \geq 1$, the $k$-{\em
  truncation} of $S$ is the
simplicial complex $T_k(S) = (V,T_k(\H))$, where $T_k(\H) = \mathcal{H} \cap P_{\leq k}(V)$.
We say that $S$ is a truncated boolean representable simplicial complex ({\em TBRSC}) if $S\; = T_k(S')$ for some $BRSC$ $S'$ and $k \geq 1$.

It is known that not every simplicial complex is a $TBRSC$ \cite[Example 8.2.6]{brsc} and not every $TBRSC$ is a $BRSC$ \cite[Example 8.2.1]{brsc}.

To understand $TBRSCs$, we need the following definition.
Given a simplicial complex $S = (V,\H)$ of rank $r$, we define
$$\varepsilon(S) = \varepsilon(\H) = 
\{ X \subseteq V \mid \forall Y \in \H \cap P_{\leq r-1}(X)\; \forall p \in V
\setminus X \hspace{.3cm} Y \cup \{ p \} \in \H\}.$$

\bl
\label{propt}
{\rm \cite[Lemma 8.2.3]{brsc}}
Let $S$ be a simplicial complex. Then:
\bi
\item[(i)] $\varepsilon(S)$ is closed under intersection;
\item[(ii)] $L(S) \subseteq \varepsilon(S)$.
\ei
\el

Thus $\varepsilon(S)$ is a Moore family and defines consequently a BRSC, denoted by $S^{\varepsilon}$.

\bt
\label{eqtr}
{\rm \cite[Theorem 8.2.5]{brsc}}
Let $S$ be a simplicial complex of rank $r$. Then the following
conditions are equivalent:
\bi
\item[(i)] $S$ is a TBRSC;
\item[(ii)] $S = T_r(S^{\varepsilon})$.
\ei
Furthermore, in this case we have $L(S^{\varepsilon}) = \varepsilon(S)$.
\et

It follows from \cite[Section 8.2]{brsc} that $S^{\varepsilon}$ is the largest $BRSC$ on $V$ whose truncation to rank $r$ is $S$.

\section{Pairwise Balanced Designs and Their Subsystems}

In \cite[Example 3.5]{tbrsc} it is shown that there are rank 3 $TBRSCs$ which are not boolean representable (unlike rank 2, see \cite[Proposition 4.1]{tbrsc}). In this section we study the class of rank 3 $TBRSCs$ in detail. We show its connection to other important combinatorial structures, the pairwise balanced designs and partial geometries. We are led directly into a connection between rank 3 $TBRSC$ and Wilson's fundamental results \cite{Wilson1, Wilson2, Wilson3}.

A { \em pairwise balanced design} ($PBD$) is given by the following data. Let $X$ be a finite set. Let $\mathcal{L}$ be a collection of subsets
$\mathcal{L} = \{B_{i} \mid i \in I\}$ of $X$ called blocks. We assume that $|B_{i}| > 1$ for all $i \in I$. Let $v$ be a non-negative integer and $K$ a set of positive integers. The pair $(X,\mathcal{L})$ is called a $K$-PBD of size $v$ if it satisfies the following conditions.

\begin{enumerate}
  \item {$|X| = v$}

  \item {$|B_{i}| \in K$ for all $i \in I$}

  \item Every pair of distinct points $x,y \in X$ is contained in a unique block $B_{i} \in \mathcal{L}$.
\end{enumerate}

Except for the cases $(\emptyset,\emptyset), (X,\emptyset), |X|=1$ and $(X,\{X\}),|X|>1$, a $PBD$ is the same thing as a 2-partition of $X$ in the sense of \cite{craporota}. This means that the blocks partition the collection of subsets of $X$ of cardinality 2. When we use the term $PBD$ we exclude these three cases in this paper. Thus, we assume that $1 \notin K$.

The following results of \cite{craporota} describe the connection of $PBDs$ to rank 3 simple matroids. We give the details for completeness.

 \bp
 \label{matroidtopbd}
 Let $M=(V,\mathcal{H})$ be a simple matroid of rank 3. Then $(V,\mathcal{L})$ is a PBD where $\mathcal{L}$ is the set of closures of two element sets of $M$.
 \ep

 \proof
In the context of matroids, it is easy to see that if $x,y$ are distinct points of $V$, then the flat generated by $x,y$ is $Cl\{x,y\}=\{x,y\} \cup \{u \in V \mid \{x,y,u\} \notin \mathcal{H}\}$ and is a proper subset of $V$. Since $M$ has rank 3, it follows that the intersection of two distinct proper flats has cardinality at most 1. Therefore every pair of distinct elements of $V$ are in a unique block and $(V,\mathcal{L})$ is a $PBD$.
 \qed

\bp
\label{pbdtomatroid}
Let $(V,\mathcal{L})$ be a $PBD$. Let $\mathcal{H} = P_{\leq 2}(V) \cup \{X \in P_{3}(V) \mid X \nsubseteqq L, \forall L\in \mathcal{L}\}$. Then $(V,\mathcal{H})$ is a simple matroid of rank 3.
\ep

\proof
We just need verify that $(V,\mathcal{H})$ satisfies the Exchange Axiom. Let $X = \{x,y\}$ be a set of size 2 and $\{u,v,w\} \in \mathcal{H}$. If $X \subset \{u,v,w\}$ then we are done, so we can assume without loss of generality that $X \cap \{u,v\} = \emptyset$. Assume that neither $X \cup \{u\}$ nor $X \cup \{v\}$ are in $\mathcal{H}$. From the definition of $\mathcal{H}$ and $(V,\mathcal{L})$ being a $PBD$, it follows that there is an $L \in \mathcal{L}$ such that
$\{x,y,u,v\} \subseteq L$. But then $w \notin L$, since $\{u,v,w\} \in \mathcal{H}$. It follows that $w \notin X$ and that $\{x,y,w\} \in \mathcal{H}$, completing the proof. \qed

A a corollary of these propositions, we see that the lattice of flats of a rank 3 matroid is constructed as follows.

\bc
Let $(V,\mathcal{L})$ be a $PBD$. Then $\P_{\leq 1}(V) \cup \L \cup \{V\}$ is closed under intersection and is the lattice of flats of the matroid from Proposition \ref{pbdtomatroid}. Every rank 3 geometric lattice is constructed in this manner.
\ec

Some issues of terminology. In \cite{Wilson1}, Wilson calls a subset $X$ of a $PBD$ $(V,\mathcal{L})$ a {\em flat} if for all distinct points $x \neq y \in X$, the unique block $\overline{xy}$ of $(V,\mathcal{L})$ containing $x,y$ is contained in $X$. This is what Crapo \cite{Crapo1} calls a 2-closed set. The term ``flats" is also an integral part of the theory of matroids, combinatorial geometry  and the theory of $BRSC$ \cite{Oxley,brsc} where they have a different meaning. To avoid confusion, we will call flats in Wilson's sense, {\em subsystems} of a $PBD$. Pairwise balanced designs are called {\em linear spaces} \cite{linspace} by combinatorial geometers. In this context subsystems are called {\em subgeometries}. We will not use this term.

The main result of this section is that flats in Wilson's sense are indeed exactly the same as flats in the sense of the theory of $BRSC$. Let $S = (V,\mathcal{H})$ be a TBRSC of rank $r$. In Theorem \ref{eqtr} we showed how to compute the largest $BRSC$ $S^{\varepsilon}$ on $V$ whose truncation to $r$ is $S$. The next theorem gives a precise connection between $\varepsilon(\mathcal{H})$ and the lattice of flats in Wilson's sense of the $PBD$ of a rank 3 matroid.

\bt
\label{tbrscwilson}
Let $(V,\mathcal{L})$ be a $PBD$ and let $M$ be the corresponding rank 3 matroid.Then $\varepsilon(M)$  is equal to the lattice of flats, $Fl((V,\mathcal{L}))$ in the sense of Wilson, of $(V,\mathcal{L})$. That is, $\varepsilon(M)$ is equal to the lattice of subsystems of $(V,\mathcal{L})$.
\et

\proof
Write $M = (V,\mathcal{H})$. Since $M$ is a matroid of rank 3, we have
$$\varepsilon(M) = \{ X \subseteq V \mid \forall Y \in \H \cap P_{\leq 2}(X)\; \forall p \in V\setminus X \hspace{.3cm} Y \cup \{ p \} \in \H\}.$$
Let $X \in \varepsilon(M)$. Then for all distinct $x,y \in X$ and $p \in V \setminus X$, $\{x,y,p\} \in \H$. By Proposition \ref{pbdtomatroid} this implies that $p$ is not in the unique block $\overline{xy}$ of $(V,\mathcal{L})$. Therefore $\overline{xy}$ is contained in $X$ and $X$ is a subsystem of $(V,\mathcal{L})$.

Conversely, assume that $X$ is a subsystem of $(V,\mathcal{L})$. Let $x \neq y \in X$ and $p \in V \setminus X$. Then $p$ is not in $\overline{xy}$ since $X$ is a subsystem. Therefore, $\{x,y,p\} \in \H$ by Theorem \ref{pbdtomatroid} and therefore $X \in \varepsilon(M)$. \qed

Despite the simplicity of the result, we see that the theory of BRSC and TBRSC are a missing link between these theories and the theory of $PBDs$. In the next section we give examples of the connection given by Theorem \ref{tbrscwilson}.

\section{Examples} \label{PBDExs}

We look at a number of examples in this section. The book  \cite{linspace} includes an Appendix containing all $PBDs$ on at most 9 points.

\be \label{CGs}

{\bf Complete Graphs} We can identify the unique 2-$PBD$ on $V$ with the complete graph on $V$. The corresponding matroid is $M=U_{|V|,3}$ whose independent sets are $P_{\leq 3}(V)$. Clearly every subset of $V$ is a subsystem in this case. Therefore, $M^{\varepsilon} = U_{|V|,|V|}$, the uniform matroid on $V$.
\ee

\be \label{NP}

{\bf Near Pencils} Let $V =\{0,1, \ldots, n\}$. Let $\L=\{\{0,i\}|i=1, \ldots, n\} \cup \{\{1, \ldots, n\}\}$. That is, $\mathcal{L}$ consists of the block $\{1, \ldots, n\}$ and all 2-sets containing 0. Then $NP(n) = (V, \L)$ is a $\{2,n\}-PBD$ called a near pencil. The corresponding matroid $M$ has as set of bases all 3-sets that contain 0. It is easy to check that the flats of $M$ are the empty set, all singletons, all the blocks of $NP(n)$ and $V$. A straightforward calculation then shows that $\varepsilon(M) = L(M)$ and $M^{\varepsilon} = M$.

\ee

\be

{\bf Projective spaces} \label{Proj}

Let $\mathbb{F}_q$ be the field of order $q$ and let $\mathbb{F}_{q}^{n+1}$ be an $n+1$ dimensional vector space over $\mathbb{F}_q$. We can consider projective $n$-dimensional space over  $F_q$ to be a PBD $P_{n,q}$ as follows. The 1-dimensional subspaces of $\mathbb{F}_{q}^{n+1}$ are the points and the 2-dimensional subspaces of $\mathbb{F}_{q}^{n+1}$ are the blocks. Incidence is given by containment. It is well known that $P_{n,q}$ is a $PBD$ on a set of size $q^{n} +q^{n-1} + \ldots + q+1$ points and $K=\{q+1\}$. The corresponding matroid $M = (V,\H)$ has basis all sets of 3 lines through the origin that are not co-planar. $\varepsilon(\H)$ is easily seen to be all the projective subspaces of $\mathbb{F}_{q}^{n+1}$ in the usual sense of projective geometry.

\ee

\be

{\bf Affine spaces} \label{Aff}

Let $V=\mathbb{F}_q^{n}$ be an $n$-dimensional space over $\mathbb{F}_q$. Affine $n$-space is the structure whose vertices are $V$ and whose blocks are all the cosets of the form
$W + v$, where $W$ is a one dimensional subspace of $V$ and $v \in V$. This is a $PBD$ with $q^n$ vertices and $K = \{q\}$. The corresponding matroid consists of all 3-sets of non-collinear points. The subsystems of the $PBD$ are the usual affine subspaces.

\ee

The above $PBDs$ all have the properties that their lattice of flats is a geometric lattice, equivalently the $BRSC$ defined by the lattice of flats of the $PBD$ is a matroid. We present examples that do not have this property. The first example was constructed by Marshall Hall in 1943 on a set of size 21. See \cite{Hallplanes}, page 236 for details.

\be

Let $V = \{1,2,3,4,5,6\}$. Let $\mathcal{L}$ consist of the sets $\{1,2,3\}, \{1,5,6\}$ and $\{3,4,5\}$ together with all two element sets not contained in any of these. This defines a $\{2,3\}$-PBD. It is easy to see that in the corresponding matroid $(V,\H)$, $H = P_{\leq 3}(V) \setminus \{ \{1,2,3\}, \{1,5,6\}, \{3,4,5\} \}$ and $\varepsilon(\H) = P_{\leq 1}(V) \cup \L \cup \{ \{2,4,6\},V\}$. Therefore,
both $\{1,3,6\}$ and $\{1,2,4,6\}$ are bases in the $BRSC$ corresponding to this $PBD$ and thus the $BRSC$ is not a pure simplicial complex and in particular, not a matroid.

\ee


\be 

The next example is related to the classical Desargues configuration. Let $K_5$ be the complete graph on 5 points. The graphic matroid $G(K_5)$ on $K_5$ has all its subforests as independent sets. We let $D$ be the truncation of $G(K_5)$ to rank 3, so that $D$ has independent sets all subforests with at most 3 edges. Since matroids are closed under truncation, $D$ is a matroid.

The corresponding $PBD$ has as points, the edges of $K_5$ and as blocks all pairs of parallel edges and all 3 sets that form a triangle. Since the latter can be
identified as the lines of the Desargues configuration, we call $D$ the Desargues matroid. By general results about truncations (\cite[Chapter 7]{Oxley}, \cite[Proposition 8.2.2]{brsc}), the lattice of flats $L(D)$ is equal to the Rees quotient of $L(G(K_{5}))$, considered as a join lattice, by the ideal of all partitions with at most two equivalence classes. Recall \cite{qtheory} that the Rees quotient of a semigroup $S$ by an ideal $I$ identifies all elements of $I$ with 0 and leaves all elements of $S \setminus I$ alone. It is not difficult to prove that $\varepsilon(D)$ is the full partition lattice and thus the corresponding $BRSC$ is $G(K_5)$.

Now we consider the ``non-Desargues" matroid. Recall \cite[Chapter 1.5]{Oxley}, that if $M$ is a matroid with set of bases $\B$ and $X$ is both a circuit and a hyperplane (that is, a flat of co-rank 1, that is, of rank one less than that of the matroid) of $M$, then $\B \cup \{X\}$ is the set of bases of a matroid called the relaxation of $M$ with respect to $X$. Any triangle in $K_5$ is indeed a hyperplane and a circuit of $D$ and fixing $T =\{34,35,45\}$ we obtain the non-Desargues matroid $N$ by relaxation of $D$ with respect to $T$. We analyze $\varepsilon(N)$ and the corresponding $BRSC$ in the next example. We note that by general facts about relaxations, the lattice of flats of $N$ is
$L(N)= P_2(T) \cup L(D) \backslash\{T\}$. Thus every flat, thought of as a subgraph of $K_5$ is a disjoint union of cliques and possibly a subset of order of 2 of $T$.

By Theorem \ref{tbrscwilson}, $\varepsilon(N)$ is the lattice of subsystems of the corresponding $PBD$. These in turn are obtained by closing subsets under the operation that for any subset $X$ of $V$ adjoins the flat of $N$ generated by any pair of distinct elements to $X$. We claim that $\varepsilon(N)$ is, by
considering a set of edges  as a subgraph of $K_5$, equal to the set of graphs on 5 points, all of whose connected components are either cliques or a 2-element subset of $T$. Clearly any such set is a subsystem. Conversely, every flat of $N$ has the required form. The flat generated by a pair of points, that is edges in $K_5$ is either that pair, if they have no point in common or they are a two element subset of $T$ or the unique triangle containing the pair if they have a point in common. By iterating this operation the required property is preserved. Thus, every subsystem has this property.

Let $E$ denote the set of edges of $K_5$. Now it is easily seen that the chain 
$$\emptyset \subset Clique(\{1,2\}) \subset Clique(\{1,2,3\}) \subset Clique(\{1,2,3,4\}) \subset E$$
is a maximal chain in $\varepsilon(N)$. But so is 
$$\emptyset \subset Clique(\{3,4\}) \subset \{34,45\} \subset Clique(\{3,4,5\}) \subset Clique(\{1,2\}) \cup Clique(\{3,4,5\}) \subset E.$$
Therefore, $\varepsilon(N)$ is not a graded lattice and in particular, not a geometric lattice and thus the $BRSC$ of $\varepsilon(N)$ is not a matroid.
\ee

\section{Wilson Monoids}

In \cite{Wilson1, Wilson2, Wilson3}, Wilson proved the existence theorem for $PBDs$ which we recall here. If $K$ is a set of positive integers, define two numbers as follows. $\alpha(K) = \text{gcd}\{k-1|k \in K\}$ and $\beta(K)=\text{gcd}\{k(k-1)|k \in K\}$. It is not difficult to prove that if $(V, \L)$ is a $K-PBD$ and $|V|=v$, then
$v-1 \equiv 0 \text{ mod }\alpha(K) \text{ and } v(v-1) \equiv 0 \text{ mod } \beta(K)$. Wilson's Theorem proves that except for a finite number of cases, if $|V|=v$
satisfies these congruential conditions, then there exists a $K-PBD$ with points $V$.

Wilson proves his theorem by combining direct constructions, that is, $PBDs$ built from finite fields, finite groups and other algebraic structures with recursive techniques to build bigger $PBDs$ from smaller pieces. Wilson implicitly defines a category of $PBDs$ by defining a notion of morphism. The self-morphisms of a $PBD \: X$ then have the structure of a monoid $W(X)$ that we call the Wilson monoid on $X$. We will explore the relationship between combinatorial properties of $X$ and
semigroup theoretic properties of $W(X)$. This leads to surprising connections between these two theories.

We begin with an example before giving formal definitions. We will call a $PBD$ {\em subsystem-free} if its only subsystems are the empty set, the singleton sets, the blocks and the whole point set. That is, a $PBD$ is subsystem-free if its only subsystems are the ones that every $PBD$ has. Equivalently, this means that the lattice of subsystems of the $PBD$ is the lattice of flats of the corresponding rank 3 matroid. Such geometries are also called non-degenerate planes \cite{Doyen1}, but we prefer the term subsystem-free. It is quite easy to see that the Fano plane is a subsystem-free $PBD$. Of course, it is a Steiner triple system (that is, a $\{3\}-PBD$).

We first build a $\{3,7\}-PBD$. We start with  the Cayley table of the group of order 7 as a Latin Square, $LS(Z_{7})$ with rows $R_{1},...,R_{7}$ and columns $C_{1},...,C_{7}$. We begin with  three blocks of size 7 consisting of the $R_{i}$, $C_{i}$ and ${i}, i=1,...,7$. We add all 49 blocks of size 3 that we obtain from $LS(Z_{7})$
of the form $\{ R_{i},C_{j},i+j(mod 7)\}$, $i,j=1,...,7$. Since any two entries of such a triple uniquely determines the third, we obtain a $\{3,7\}-PBD, PBD(Z_{7})$. A short calculation will show that this is a subsystem-free $PBD$.

We now use the technique \cite{Wilson1} to build a Steiner triple system,  on the 21 points of $PBD(Z_{7})$. We replace or ``break up" each of the three blocks of size 7 with disjoint copies of the Fano plane. It is easy to see that this is indeed a triple system $X$ on 21 points. Furthermore, the blocks that were of size 7 in $PBD(Z_{7})$ are now flats of size 7 in $X$ and thus $X$ is not a subsystem-free $PBD$.

Clearly we can use any Latin Square $L$ on $n$ points in place of $LS(Z_{7})$ and any Steiner triple system of order $n$ to build a $\{3,n\}-PBD$, $PBD(L)$ on 3$n$ points. Steiner triple systems built this way are called systems of Wilson-type in \cite{sts1921}. We will look in more detail at Wilson monoids of Steiner triple systems later in the paper.

We now define the morphisms in Wilson's sense between PBDs. We first need a non-conventional definition of inverse image of partial functions. Let $f:S\rightarrow T$ be a partial function between sets $S$ and $T$. If $S_0$ is the domain of $f$ we call $S \setminus S_0$ the {\em co-domain} of $f$. Wilson \cite{Wilson1} calls this the {\em kernel} of $f$, but we use this term for the partition on $Dom(f)$ that identifies two elements if they have the same image under $f$. We let $f_{0}:S_{0} \rightarrow T$ be the total function defined by $f$. If $A \subseteq S$, then we let $f(A)=f_{0}(A \cap S_{0})$ and if $B \subseteq T$, we define $f^{-w}(B)=f_{0}^{-1}(B) \cup (S \setminus S_{0})$. We use the notation $f^{-w}$ to denote the inverse image in the sense of Wilson and the usual notation $f^{-1}$ for the standard notion of inverse image of a partial function. Thus, the co-domain of $f$ is contained in the Wilson inverse image of any subset $B$ of $T$.

Let $X=(S,\mathcal{L})$ and $Y=(T,\mathcal{L'})$ be $PBDs$. A partial function $f:S \rightarrow T$ is called a morphism between $X$ and $Y$ if $f^{-w}(F)$ is a subsystem of $X$ for every subsystem $F$ of $Y$. Notice that by the definition of Wilson inverse image, $f^{-w}(\emptyset)$ is the co-domain of $f$ and since the empty set is a subsystem, the co-domain of any morphism is a subsystem of $X$. We define an {\em open set} of a $PBD$ to be the complement of a subsystem of $X$ and it follows that the domain of a morphism is an open set. The following straightforward lemma allows us to use the following equivalent definition of morphism in terms of open sets and the usual definition of inverse image in the rest of the paper. This also allows us to use the results of \cite{bibdtrans, bibdcont, projcont} to understand the monoid of morphisms on a $PBD$.

\bp \label{opendef}
Let $X=(S,\mathcal{L})$ and $Y=(T,\mathcal{L'})$ be $PBDs$. A partial function $f:S \rightarrow T$ is a morphism if and only if the domain of $f$ is an open subset and $f^{-1}(O)$ is an open set of $X$ for every open set $O$ of $Y$.
\ep

\proof
Let $f:S \rightarrow T$ be a morphism between $X$ and $Y$ and let $O$ be an open subset of $T$. By definition, $f^{-w}(T \setminus O)=f_{0}^{-1}(T \setminus O) \cup \text{co-domain}(f)$ is a subsystem of $X$. But $S$ is the disjoint union of $f^{-1}(O)$ and $f_{0}^{-1}(T \setminus O) \cup \text{co-domain}(f)$ and thus $f^{-1}(O)$ is an open subset of $X$.

Conversely assume that the domain of $f:S \rightarrow T$ is open and that $f^{-1}(O)$ is an open set of $X$ for every open set $O$ of $Y$. Then $f^{-w}(\emptyset)=\text{co-domain}(f)$ is a subsystem of $X$.

Now let $F$ be a non-empty subsystem of $Y$. Then $f^{-1}(T\setminus F)$ is an open subset of $S$. Clearly, the complement of $f^{-1}(T\setminus F)$ in $S$ is
$f^{-1}(F) \cup \text{co-domain}(f)=f^{-w}(F)$. Therefore, $f^{-w}(F)$ is a subsystem of $X$ and thus $f$ is a morphism. \qed

\bc
Let $f:S \rightarrow T$ and $g:T \rightarrow U$ be morphsims of $PBDs$. Then $gf:S \rightarrow U$ is a morphism of $PBDs$.
\ec

\proof
Since $\text{Dom}(g)$ is an open set and $f$ is a morphism if follows that
$\text{Dom}(gf)=f^{-1}(\text{Dom}(g))$ is an open set. Also, if $O$ is an open subset of U, then $(gf)^{-1}(O)=f^{-1}(g^{-1}(O))$ is an open set of $S$ since both $f$ and $g$ are morphisms. It follows from Proposition \ref{opendef} that $gf:S \rightarrow U$ is a morphism. \qed

This allows us to define $\mathcal{PBD}$ to be the category whose objects are $PBDs$ and whose morphisms are those defined in this section. In particular, for every $PBD$ $X=(S,\mathcal{L})$, we define its Wilson monoid $W(X)$ to be the monoid of all morphisms from $X$ to itself. We will elucidate the structure of $W(X)$ in the rest of this section.

We interpret Proposition \ref{opendef} as follows. The collection $\mathcal{O}$ of open subsets of a $PBD$ is closed under unions and contains the co-points, that is, the sets of cardinality one less than $V$, as well as the empty set and the whole set. Thus $\mathcal{O}$ satisfies all the axioms of a topology except possibly closure under intersection. In this ``generalized topology", Proposition \ref{opendef} says that the Wilson morphisms between $PBDs$ are precisely the partial continuous functions. It was this analogy that lead Dinitz and Margolis to call such partial functions between arbitrary incidence structures continuous partial functions \cite{bibdcont,bibdtrans}. We view the category $\mathcal{PBD}$ as a natural generalization of the category of topological spaces.

We begin with the following very important proposition of Wilson \cite[Proposition 7.1]{Wilson1} that gives a characterization of morphisms by their effect on direct image on {\em blocks} of a $PBD$. Thus Wilson self-morphisms are a special kind of endomorphism of a $PBD$. We give the proof for purposes of completeness.

\bp \label{morph1}
Let $X=(S,\mathcal{L})$ and $Y=(T,\mathcal{L'})$ be $PBDs$. A partial map $f:S \rightarrow T$ is a morphism if and only if the domain of $f$ is open and for every block
$B \in \mathcal{L}$, either (i) $|f(B)| \leq 1$ or (ii) $f$ is defined on all of $B$, is one-to-one on $B$ and there is a (necessarily unique) block $B' \in \mathcal{L'}$ such that $f(B) \subseteq B'$.
\ep

\proof
Assume that the domain of $f$ is open and satisfies conditions (i) or (ii) for every block $B \in \mathcal{L}$. Let $F$ be a subsystem of $Y$ and let $E=f^{-w}(F)$. If $|E| \leq 1$, then $E$ is a subsystem of $X$. Assume then, that $x_{1},x_{2}$ are two distinct points of $E$ and let $B$ be the unique block of $X$ that contains $x_{1},x_{2}$. If (i) holds, then either $f(B)$ is the empty set or $f(B)=y$ for some $y \in F$. In both cases, $B \subseteq f^{-w}(F)=E$. If $|f(B)|>1$, then by (ii), $f$ is defined on all of $B$ and is one-to-one on $B$ and there is a block $B'$ of $Y$ such that $f(B) \subseteq B'$. $B'$ contains the two distinct points $f(x_{1}),f(x_{2})$ of $f(E)=F$. Since $F$ is a subsystem of $Y$, $B' \subseteq F$ and thus $B \subseteq f^{-1}(B') \subseteq f^{-w}(F)=E$. Therefore, $f$ is a morphism.

Conversely, assume that $f:S \rightarrow T$ is a morphism. Then the domain of $f$ is open. Let $B$ be a block of $X$. Since block sizes are greater than 1, if $f$ is either
not defined on all of $B$ or is not one-to-one on $B$, then there are two distinct points $x_{1},x_{2}$ in $B$ such that $|\{f(x_{1}),f(x_{2})\}| \leq 1$. Therefore the
 set $F=\{f(x_{1}),f(x_{2})\}$ is a subsystem of $Y$ and since $f$ is a morphism, $f^{-w}(F)$ is a subsystem that contains the two distinct points $x_{1},x_{2}$. Thus,
  $B \subseteq F$ and thus $f(B) \subseteq \{f(x_{1}),f(x_{2})\}$ and it follows that (i) holds. It follows that if (i) doesn't hold then $f$ is defined on all of
$B$ and is one-to-one on $B$. Therefore, for two distinct points, $x_{1},x_{2}$ in $B$, $|\{f(x_{1}),f(x_{2})\}|=2$ and thus lie in a unique block $B'$ of $Y$. Since $f$ is a morphism, $x_{1},x_{2}$ are contained in the subsystem $f^{-w}(B')$ and it follows that $B \subseteq f^{-w}(B')$. Since $f$ is defined on all of $B$,
$B \subseteq f^{-1}(B')$. Therefore, $f(B) \subseteq B'$ and (ii) holds. \qed

A morphism between $X=(S,\mathcal{L})$ and $Y=(T,\mathcal{L'})$  is called an {\em open morphism} if and only if the image of every subsystem of $X$ is a subsystem of $Y$. A proof similar to that of Proposition \ref{morph1} proves the next proposition.

\bp \label{morph2}
Let $X=(S,\mathcal{L})$ and $Y=(T,\mathcal{L'})$ be $PBDs$. A morphism $f:S \rightarrow T$ is an open morphism if and only if for every block $B \in \mathcal{L}$, either $|f(B)| \leq 1$ or $f(B)$ is a block of  $Y$.
\ep

Proposition \ref{morph1} and Proposition \ref{morph2} show that the blocks of a $PBD \, X =(S,\mathcal{L})$ form a weakly-preserved cover of the action of $W(X)$. That is, the union of all the blocks is all the points, and the image of any block under the action of any element of $W(X)$ is contained in (a not necessarily unique, if case (i) of Proposition \ref{morph1} holds). Weakly-preserved covers play an important part in Zieger's proof of the Krohn-Rhodes Theorem \cite{Zeiger1}. In this and a future paper, we exploit the properties in these two propositions and use the interaction of the combinatorics of $X$ and the geometry of the actions of $W(X)$ on points, open sets and blocks to study various decompositions: one and two-sided wreath products, triangular products \cite{qtheory} of $W(X)$ and its semiring of subsets $P(W(X))$.

The following corollary follows easily from Proposition \ref{morph1} and Proposition \ref{morph2}. An automorphism of a $PBD$ is a permutation on the points that sends blocks to blocks. A partial constant function $f:S \rightarrow S$ is a partial function such that $|f(S)| \leq 1$.

\bc \label{auto}

Let $X=(S,\mathcal{L})$ be a PBD.
\begin{itemize}
  \item[(i)] {A permutation $f:S \rightarrow S$ is a morphism if and only if $f$ is an automorphism of $X$.}
  \item[(ii)] {A partial constant map $f:S \rightarrow S$ is a morphism if and only if f is the empty function or the domain of $f$ is a non-empty open subset $O$ of $S$ and its image is a point $p \in S$.}
\end{itemize}
\ec

Let $X=(S,\mathcal{L})$ be a $PBD$. If the morphism $f:S \rightarrow S  \in W(X)$ is a non-empty partial constant function, then we write $f=(p,O)$ if the domain of $f$ is the non-empty open subset $O$ of $S$ and its image is $\{p\}$. We write $\theta$ for the empty function. Clearly, the collection of all partial constant functions in $W(X)$ is an ideal in $W(X)$. We identify this ideal as the unique 0-minimal ideal of $W(X)$ and compute its structure as a 0-simple semigroup and how it sits inside $W(X)$ as an ideal.

Let $(p',O'),(p,O)$ be two partial constant functions in $W(X)$. Clearly, $${(p',O')(p,O) = (p',O) \text{ if } p \in O'}$$ and the empty function otherwise. If we define a pairing $<,>: \mathcal{O} \times S \rightarrow \{0,1\}$, where $\mathcal{O}$ is the collection of non-empty open sets of $X$ by $<O',p>\; =1 \text{ if } p \in O'$ and 0 otherwise, then we can write the product above as $(p',O')(p,O) = (p'<p,O'>,O)$,where we identify multiplication by 0 as giving the empty function. It is straightforward then to see that the ideal $I(X)$ of partial constant maps of $W(X)$ is isomorphic to the Rees matrix semigroup \cite[A4]{qtheory} $\mathcal{M}^{0}(\{1\}, S,\mathcal{O},<,>)$.

We now note that the natural left action of $W(X)$ on $S$, that is, $fp =f(p), f \in W(X), p \in S$ by partial functions has an ``adjoint" right action on $\mathcal{O}$. Namely, let $f \in W(X)$ and define the partial function $\bar{f}:\mathcal{O} \rightarrow \mathcal{O}$ acting on the right of $\mathcal{O}$ by $O\bar{f}=f^{-1}(O)$ if this set is non-empty and undefined otherwise. Adjointness means that for the pairing $<,>$ defined above, we have ${<O \bar{f},p>\; =\; <O,fp>}$ for all $O \in \mathcal{O}, p \in S$. In the language of semigroup theory, this means that $W(X)$ is the translational hull of the 0-simple semigroup $I(X)$. See \cite[Section 5.5]{qtheory} for a general introduction to the translational hull of a finite 0-simple semigroup.

We note that the pairing $<,>$ (also known as the structure matrix of the 0-simple semigroup) is
{\em reduced}. This means that for all distinct $p,q \in S$ there is an $O \in \mathcal{O}$ such that $<p,O>\; \neq\; <q,O>$ (since any two points are in exactly one block and we are assuming that there are at least two blocks) and for each $O \neq O' \in \mathcal{O}$, there is a $p \in S$ such that $<p,O>\; \neq\; <p,O'>$. If we think of $<,>$ as a $|S| \times |\mathcal{O}|$ matrix over $\{0,1\}$, then reduced means that distinct rows (columns) are not equal to one another. 0-simple semigroups over the trivial group with reduced structure matrices are precisely the congruence-free 0-simple semigroups  and along with the semigroups of order 2 and finite simple groups, form the class of all finite congruence-free semigroups \cite[Theorem 4.7.17]{qtheory}. A semigroup $S$ is called Generalized Group Mapping (GGM) if it has a unique 0-minimal ideal $I(S)$ which is a 0-simple semigroup and such that $S$ acts faithfully on both the left and right of $I(S)$ by left and right multiplication \cite[Chapter4]{qtheory}. More precisely, $S$ acts faithfully by partial functions on any $\mathcal{L}$ and $\mathcal{R}$ class in $I(S) \setminus \{0\}$, which means both the left and right Sch$\rm \ddot{u}$tzenberger representations on the $\mathcal{J}$-class $I(S) \setminus \{0\}$ are faithful. An important theorem says that if the maximal subgroup of $I(S) \setminus \{0\}$ is trivial, then $S$ is GGM if and only if $I(S)$ is a congruence-free 0-simple semigroup and $S$ is a subsemigroup of the translational hull of $I(S)$ \cite[Sections 4.6, 5.5]{qtheory}. We summarize all of this discussion in the following Theorem. By ``non-trivial $PBD$" we mean one that contains at least two blocks, that is, it is not the $PBD$ $(S,\{S\})$.

\bt
Let $X=(S,\mathcal{L})$ be a non-trivial PBD and $W(X)$ its Wilson monoid of continuous functions. Let $\mathcal{O}$ be the collection of non-empty open subsets of $X$. Then $W(X)$ is a GGM semigroup with unique 0-minimal ideal $I(X)$ isomorphic to the congruence-free Rees matrix semigroup $\mathcal{M}^{0}(\{1\}, S,\mathcal{O},<,>)$, where $<O,p>$ is 1 if $p \in O$ and 0 otherwise. Furthermore, $W(X)$ is isomorphic to the translational hull of $I(X)$.

\et

\subsection{Examples}

In this subsection, we look at the Wilson monoids of the examples in section \ref{PBDExs}.

\be

{\bf Complete Graphs}

\ee
We saw in Example \ref{CGs} that the complete graph on a set $V$ is the unique {2}-$PBD$ on $V$. We saw that every subset of $V$ is both a subsystem and hence every subset is also open. Therefore every partial function is a morphism and the Wilson monoid of the complete graph is the monoid of all partial functions.

The next two examples show that as one might expect, projective spaces and affine spaces have many continuous maps arising from the ambient monoid of matrices.

\be
{\bf Projective Spaces}\label{ProjWil}
\ee
As in Example \ref{Proj} we consider Projective $n$-dimensional space over  $\mathbb{F}_q$, the field of order $q$ to be the PBD $P_{n,q}$ whose points are the
1-dimensional subspaces of $\mathbb{F}_{q}^{n+1}$ and the 2-dimensional subspaces of $\mathbb{F}_{q}^{n+1}$ are the blocks.
 A semilinear function $f:\mathbb{F}_{q}^{n+1} \rightarrow \mathbb{F}_q^{n+1}$ is a function such that $f(v+w) =f(v)+f(w)$ for
 all $v,w \in \mathbb{F}_q^{n+1}$ and $f(cv)=\sigma(c)f(v)$ for all $c \in \mathbb{F}_{q}, v \in \mathbb{F}_q^{n+1}$ and where
$\sigma:\mathbb{F}_{q} \rightarrow  \mathbb{F}_q$ is a fixed automorphism of $\mathbb{F}_q$. An invertible semilinear map clearly sends points of $P_{n,q}$ to itself and
 preserves incidence, so defines an
  automorphism (also called a collineation) of $P_{n,q}$. The fundamental theorem of projective geometry \cite{geoalg} states conversely, that every automorphism of
   $P_{n,q}$ is induced by an invertible semilinear map. It follows from Corollary \ref{auto}, that the group of units of $W(P_{n,q})$ is this same group.

More generally, let $f:\mathbb{F}_q^{n+1} \rightarrow \mathbb{F}_q^{n+1}$ be an arbitrary semilinear map. The kernel $ker(f)$ of $f$, that is, the set of all $v \in \mathbb{F}_q^{n+1}$ sent to 0 by $f$ is a subspace of $\mathbb{F}_q^{n+1}$. Let $[ker(f)]$ denote the subspace of $P_{n,q}$ associated to $ker(f)$. We now define
$\bar{f}:P_{n,q} \rightarrow P_{n,q}$ a partial function with domain $P_{n,q} \setminus [ker(f)]$. That is, the domain of $\bar{f}$ consists of the one dimensional
subspaces of $\mathbb{F}_{q}^{n+1}$ not contained in $ker(f)$. Therefore, if $v$ is such a line, $f(v)$ is also a one dimensional subspace of $\mathbb{F}_{q}^{n+1}$ and we define $\bar{f}(v)$ to be the point $f(v)$ of $P _{n,q}$.

We note that the domain of $\bar{f}$ is an open subset of $P_{n,q}$, being the complement of a subspace of $P_{n,q}$. Now let $b$ be a block of $P_{n,q}$, that is a 2-dimensional subspace of $\mathbb{F}_{q}^{n+1}$. If $b \subseteq ker(f)$, then $\bar{f}(b)$ is the empty set. If the intersection of $b$ and $ker(f)$ is one dimensional, then $f(b)$ is a one dimensional subspace of $\mathbb{F}_{q}^{n+1}$, so that $\bar{f}(b)$ is a point of $P_{n,q}$. Finally, if $b \cap ker(f)=\{0\}$, then $f$ maps $b$ one-to-one onto the 2 dimensional space $f(b)$ and induces a bijection on the one dimensional subspaces from those of $b$ to those of $f(b)$. Therefore, in this case, $\bar{f}$ is one-to-one on the points of $b$ considered as a block in $P_{n,q}$. It follows from Proposition \ref{morph1} that $\bar{f}$ is an element of $W(P_{n,q})$.

See \cite{projcont} where it is proved that the monoid of continuous functions on a design defined on projective space is the monoid of all projective matrices over the corresponding field.

\be

{\bf Affine Spaces}\label{AffWil}

\ee

In Example \ref{Aff} we defined $n$-dimensional affine space $A(n,q)$ over the field $\mathbb{F}_q$ to be the $PBD$ whose points are the elements of $\mathbb{F}_{q}^{n}$ and whose blocks are all the cosets of one-dimensional spaces of $\mathbb{F}_{q}^{n}$. Let $M$ be an $n \times n$ matrix over $F_q$. $M$ acts on $\mathbb{F}_{q}^{n}$ and if $l$ is a one-dimensional subspace of $\mathbb{F}_{q}^{n}$ and $a \in \mathbb{F}_{q}^{n}$, then $M(l+a)=Ml+Ma$. Since the latter is either a point or is a block which is a bijective image of $l+a$, $M$ defines a total continuous function by Proposition \ref{morph1}. More generally, any affine function on $\mathbb{F}_{q}^{n}$, that is a function $f:\mathbb{F}_{q}^{n} \rightarrow \mathbb{F}_{q}^{n}$ of the form $f(v)=Mv + w$, where $M$ is an $n \times n$ matrix over $\mathbb{F}_q$ and $w$ is a fixed element of  $\mathbb{F}_{q}^{n}$ defines a continuous function on $A(n,q)$. We leave the problem of determining the full monoid $W(A(n,q))$ for later work.

\section{Group Divisible Designs and PBDs of Split Wilson Type}

Morphisms between $PBDs$ are important in that they allow a very general scheme to build large designs from smaller ones. This plays a crucial role in Wilson's proof that the easy congruential necessary conditions for the existence of designs are eventually sufficient.

Let $X=(S,\mathcal{L})$ and $Y=(T,\mathcal{L}')$ be $PBDs$ and let $f:S \rightarrow T$ be a morphism. The key is that for $B$ a block of $Y$, $f^{-1}(B)$ is either empty or a group divisible design ($GDD$), a concept that we now recall.

A $GDD$ is a triple $X=(S,\mathcal{G},\mathcal{L})$, where $S$ is a finite set, $\mathcal{G}$ is a partition of $S$ and $\mathcal{L}$ is a set of subsets of $S$ of size at least 2. Elements of $\mathcal{G}$ are called {\em groups} and elements of $\mathcal{L}$ are called {\em blocks}. It is required that every distinct pair of points $x,y \in S$, is contained in either a unique group or a unique block, but not both. If $\mathcal{G}'$ is the set of groups of size at least 2, then
$(S,\mathcal{L} \cup \mathcal{G}')$ is a $PBD$. Conversely, if $(S,\mathcal{L})$ is a $PBD$ and $\mathcal{G}'$ is a collection of blocks of $\mathcal{G}'$ that is a partial partition of $S$ (that is, a collection of non-empty disjoint subsets of $S$), then $(S,\mathcal{G},\mathcal{L})$ is a $GDD$, where $\mathcal{G}$ is the partition of $S$ consisting of the elements of
$\mathcal{G}'$ together with all the singleton subsets of elements of $S$ not in the union of the elements of $\mathcal{G}$. These operations are clearly inverses and thus a $GDD$ is the same thing as a $PBD$ with a distinguished partial partition of $S$. A subsystem of a $GDD$ is a subsystem of its corresponding $PBD$.

A $GDD$ is {\em uniform} if all its blocks have the same size. A {\em transversal} of a $GDD$ is a block $Y$ that meets every group in precisely one point. That is, a transversal is a system of distinct representatives for the groups of the $GDD$. A $(k,m)$-transveral design, $TD(k,m)$ is a uniform $GDD$ in which all blocks have size $k$ and there are $k$ groups each with $m$ elements. Thus a $TD(k,m)$ has $km$ points and each block is a transversal. Conversely, if $X$ is a $GDD$ with at least 3 groups, such that every block is a transversal, then $X$ is a $TD(k,m)$ for some $m$. \cite[Theorem 6.2]{Wilson1}.

\be

Let $L$ be a Latin square of order $m$, that is an $m \times m$ matrix $L$ with entries in $\{1,...,m\}$ such that each entry appears precisely once in every row and column of $L$. Let $S=\{R_{1},\ldots,R_{m}\} \cup \{C_{1},\ldots,C_{m}\} \cup \{1,...,m\}$, be a set of size $3m$. We let $\mathcal{G}$ be the partition of $S$ into these three sets of size $m$ and we let $\mathcal{L}=\{\{R_{i},C_{j},L(i,j)\}|1 \leq i,j \leq m\}$. Since $L$ is a Latin square, any two elements of a triple in $\mathcal{L}$ uniquely determine the third element and thus $(S,\mathcal{G},\mathcal{L})$ is a $TD(3,m)$ It is known that every $TD(3,m)$ is constructed this way.

\ee

The following is part of Theorem 8.1 of \cite{Wilson1}. We give the proof for purposes of completeness.

\bl \label{gdd}

Let $X$ and $Y$ be $PBDs$ on the sets $S$ and $T$ respectively and let $f:S \rightarrow T$ be a morphism. Then the co-domain $D$, that is the set of points on which $f$ is not defined, is a subsystem of $X$ as is $D \cup f^{-1}(y)$ for all $y \in Y$. Let $B$ be a block of $Y$ such that $f^{-1}(B)$ is not empty. Let $Z = f^{-1}(B)$,
$\mathcal{G}=\{f^{-1}(y)\mid y \in (B \cap f(S))\}$ and let $\mathcal{L}$ be the set of blocks of $X$ that are contained in $Z$ and that intersect every class of $\mathcal{G}$ in at most one point. Then $(Z,\mathcal{G},\mathcal{L})$ is a GDD.

\el

\proof Since $f$ is a morphism $Dom(f)$ is an open subset of $S$. Therefore the co-domain, $D=S \setminus Dom(f)$ is a subsystem of $X$. Let $y \in Y$. Then $\{y\}$ is a subsystem of $Y$ and thus $f^{-w}(y)=D \cup f^{-1}(y)$ is a subsystem of $X$.

Let $x_{1},x_{2}$ be two distinct points of $Z$. By definition they can not both be in  some group in $\mathcal{G}$ and a block in $\mathcal{L}$. Furthermore, since $X$ is a $PBD$, $x_{1},x_{2}$ can be in at most one block in $\mathcal{L}$. We claim that if $x_{1},x_{2}$ are in two different groups $f^{-1}(y_{1}),f^{-1}(y_{2})$ of $\mathcal{G}$, then there is some block $b \in \mathcal{L}$ containing them.

Let $b$ be the unique block of $X$ containing $x_{1},x_{2}$. Since $B$ is a block of $Y$ and $f$ is a morphism, $f^{-w}(B)=D \cup Z$ is a subsystem of $X$ containing $x_{1},x_{2}$ and thus, $b \subseteq D \cup Z$. If $b$ contained a point $x_3$ of $D$, then $b$ would be contained in the subsystem $f^{-w}(y_{1})= D \cup f^{-1}(y_{1})$ since $x_{1},x_{3}$ are in this subsystem. This implies that $x_2 \in  D \cup f^{-1}(y_{1})$ and since $x_2$ is in the domain of $f$, it must be in the group $f^{-1}(y_{1})$ contradicting the assumption that $x_{1},x_{2}$ are in two distinct groups. Therefore, $b \subseteq Z$. If $b$ contained two points in the same group $f^{-1}(y)$ of $\mathcal{G}$, then $b$ would be contained in the subsystem $f^{-w}(y)=D \cup f^{-1}(y)$ again contradicting that $x_{1},x_{2}$ are in different groups. Therefore, $b \in \mathcal{L}$ and $(Z,\mathcal{G},\mathcal{L})$ is a $GDD$. \qed

The converse of Lemma \ref{gdd} is also true. That is, if $Y$ is a $PBD$ and there is a collection of suitable sized $GDDs$, one for each block of $Y$, and $PBDs$ that play the role of $D$ and $f^{-1}(y)$ in the above proof, then there exists a $PBD$ $X$ and a morphism $f:X \rightarrow Y$ that respects this data as in Lemma \ref{gdd}. This allows one to build a $PBD$ $X$ from a $PBD$ $Y$ and a collection of suitable $GDDs$, glued together by a morphism from $X$ to $Y$. See \cite[Theorems 8.1, 8.2]{Wilson1} for details. These results are among the most important ways to build large collections of designs and show why morphisms are an important part of the theory of $PBDs$.

We now study idempotents and regular elements in Wilson monoids. Recall that a regular element $s$ of a semigroup $S$ is an element such that there exists $t \in S$ such that $sts=s$.

\bl \label{idempotent}

Let $X$ be a PBD and let $e$ be an idempotent in $W(X)$. Then the image of $e$ is a subsystem of $X$.

\el

\proof Let $x_{1} \text{ and }x_{2}$ be two distinct elements of the image of $e$ and let $b$ be the (unique) block of $X$ containing these points. Since $e$ is an idempotent, it follows that $e(x_{1})=x_{1} \text{ and } e(x_{2})=x_{2}$.  It follows from Proposition \ref{morph1}, that $e$ is defined on all of $b$ and is one-to-one on $b$ and there is a block $b'$ of $X$ such that $e(b) \subseteq b'$. As $x_{1} \text{ and }x_{2}$ are both in $e(b)$, it follows that $b'=b$ (by uniqueness) and thus $e(b)=b$. Therefore $b$ is contained in the image of $e$, which is therefore a subsystem. \qed

\bc \label{regimage}

Let $X$ be a PBD and let $f$ be a regular element in $W(X)$. Then the image of $f$ is a subsystem of $X$.

\ec

\proof Let $g \in W(X)$ be such that $fgf=f$. Then $e=fg$ is an idempotent in $W(X)$ and it is easy to prove that the image of $f$ is equal to the image of $e$. The result follows from Lemma \ref{idempotent} that the image of $f$ is a subsystem of $X$. \qed

The following example shows that despite having proved that the range of an idempotent morphism is a subsystem, it need not be an open map. That is, it need not send every subsystem onto a subsystem.

\be

Consider the $PBD$ $(\{0,1,2,3\},\{\{1,2,3\},\{0,1\},\{0,2\},\{0,3\}\})$. Then the map $f:\{0,1,2,3\}\rightarrow \{0,1,2,3\}$ defined by $f(0)=f(1)=1,f(2)=2,f(3)=3$ is an idempotent morphism by Proposition \ref{morph1}, but is not an open morphism since the image of the subsystem $\{0,2\}$ is not a subsystem.

\ee


%

Certain idempotent morphisms that we call {\em split idempotents} allow us to split a $PBD$ over a proper subsystem in the sense we now describe. Let $X$ and $Y$ be $PBDs$ on the sets $S$ and $T$ respectively and let $f:S \rightarrow T$ be a surjective morphism. A {\em section} of $f$ is a subsystem $F$ such that $f|_{F}$ is a bijection. If $f$ is an open morphism, then $f|_{F}$ is an isomorphism. In this case, if $g:Y \rightarrow X$ is the inverse morphism of $f|_{F}$, then $ e = gf:X \rightarrow X$ is an open idempotent in $W(X)$ with image $F$, which we identify with $Y$. Clearly, $e$ is an open self-morphism. In general, we call an open idempotent $e =e^{2} \in W(X)$ a split idempotent. This leads to the following definition.

\bd

A PBD X is called split of Wilson type, if W(X) contains a split idempotent e with range a subsystem Y with $1 < |Y| < |X|$. We usually just write ``$PBD$ of Wilson type."

\ed

Recall that a small monoid is a monoid that is the disjoint union of its group of units and a unique 0-minimal ideal that is a 0-simple semigroup. As we have seen that the set of all partial constant maps of a Wilson monoid form the unique 0-minimal ideal and is a 0-simple semigroup, it follows that if a $PBD$ is of Wilson type then its Wilson monoid is not small.

\section{Steiner Triple Systems of Order up to 19}

$PBDs$ of Wilson type are on the one hand rare among all $PBDs$ but are powerful enough to construct a wide range of $PBDs$ and be counted efficiently \cite{Wilson1, Wilson2, Wilson3}. A Steiner triple system ($STS$) is a $PBD$ with all blocks of size 3. That is an $STS$ is a $(v,3,1)$ Balanced Incomplete Block Design ($BIBD$). It is well known that an $STS$ exists if and only if $v\equiv 1 \text{ (mod }6) \text{ or } v\equiv 3\text{ (mod }6)$. In  this subsection we survey $STS$ of size up to 19 and their Wilson monoids.

$\bf v=3$

If $v=3$, then the unique $STS, I_3$ up to isomorphism is the trivial $STS$ with the set of points as the unique block. The open sets are the empty set, the sets of cardinality 2 and the whole set. By Proposition \ref{morph1}, $W(I_{3})$ consists of the symmetric group $S_3$ as group of units, all the maps $(p,O)$ of rank 1, where $p$ is a point and $O$ is a non-empty open set (see Corollary \ref{auto}) and the empty function. It is straightforward to compute that there are 19 elements in $W(I_{3})$.

$W(I_{3})$ is s small monoid. That is, it consists of a group of units and a unique 0-minimal ideal which is a 0-simple semigroup. We will shortly see that generically the Wilson monoid of an $STS$ is a small monoid.

$\bf v=7$

It is well known that the unique $STS$ up to isomorphism on 7 points is the Fano plane, which is isomorphic to the projective plane $P_{2,2}$ over the field of order 2. As a $PBD$, the Fano plane is a $(7,3,1)$-$BIBD$. We can identify its point set $V$ with the seven non-zero elements of $\mathbb{F}_{2}^3$.The subsystems of $P_{2,2}$ are the empty set, the points, the seven lines and the whole point set. The open sets are the complements of these.

As for any $BIBD$, every Wilson self-map on the Fano plane is open. Thus the possible ranges of Wilson maps, are the subsystems. The group of units of $W(P_{2,2})$ is the collineation group of $P_{2,2}$ which is well known to be the simple group $PSL(3,2)$, the projective special linear group of order 168. There are 15 open sets and thus the unique 0-minimal ideal of $W(P_{2,2})$ has order 106 = (15x7)+1. It follows from the description in Example \ref{ProjWil} that every linear transformation $M$ on $\mathbb{F}_{2}^3$ restricts to a Wilson map $f_{M}:V \rightarrow V$ with domain $V - Ker(M)$. If the rank of $M$ is 2, then the image of $f_M$ is a block. We will now show that every Wilson map with image a block is of this form.

Assume that $f:V \rightarrow V$ is a Wilson map with range a block $b$ of the Fano plane, so $b$ consists of the non-zero elements of a 2-dimensional subspace of $F_{2}^{3}$.
By Propositions \ref{morph1} and \ref{morph2} and sections 7-9 of \cite{Wilson1}, the domain of $f$ has 6 points and $f^{-1}(b)$ is a transversal design $TD(3,2)$. This means that for every $v \in b$, then $f^{-1}(v)$ has two points and that $f^{-w}(v)=f^{-1}(v) \cup \{p\}$ is a block of the Fano plane, where $p$ is the unique point not in the domain of $f$. Thus, $f^{-w}(b)$ is the pencil on $p$, that is, the 3 blocks of the Fano plane that pass through $p$. By elementary linear algebra, there is a linear transformation $M:\mathbb{F}_{2}^{3}\rightarrow \mathbb{F}_{2}^{3}$ with kernel $\{0,p\}$, range $b \cup \{0\}$ and such that the inverse image of points of $b$ are the non-zero cosets of $\{0,p\}$. It follows that $f=f_{M}$.

It is well known that the linear transformations of rank 2 of $\mathbb{F}_{2}^{3}$ form a $\mathcal{J}$-class of the monoid of all linear transformations of $\mathbb{F}_{2}^{3}$. The maximal subgroup of this $\mathcal{J}$-class is the general linear group $Gl_{2}(2)$, which is isomorphic to the symmetric group on 3 points. Since there are seven subspaces of $\mathbb{F}_{2}^{3}$ of dimension 1 and 2, and each pair can serve as the kernel and range of a linear transformation, there are $7\times 7 \times 6=294$ linear transformations of rank 2 over
$\mathbb{F}_{2}^{3}$. It follows together with the count above of the group of units and elements of rank at most 1 in $W(P_{2,2})$ that $|W(P_{2,2})|=568$.

$\bf v=9$

It is known that up to isomorphism the unique $STS$ on 9 points is the affine plane $AG(2,3)$ over the field of order 3. The subsystems are all affine subspaces of $\mathbb{F}_{3}^2$ including the empty set and the open sets are their complements. As mentioned in Example \ref{AffWil}, every affine function on $\mathbb{F}_{3}^2$ defines a Wilson map on $AG(2,3)$. An argument similar to the one in the previous example shows that $W(AG(2,3))$ consists of the affine functions together with all the partial constant maps.

Before continuing, we need two results. The first is a well known result about BIBDs generalizing Fisher's inequality. See Proposition 4.1 of \cite{Wilson1}, for a proof.

\bp \label{subbound}

Let X be a (v,k,1)-BIBD, $k\geq 2$. If X has a subsystem of order $u < v$, then $v \geq (k-1)u+1$. In particular, if X is an STS, then $v \geq 2u+1$.

\ep

Let $X$ be a $BIBD$. By Propositions \ref{morph1} and \ref{morph2}, every element of $W(X)$ is an open morphism. The following additional property follows immediately from Proposition 9.2 and Theorem 9.3 of \cite{Wilson1}.

\bt \label{uniform}

Let $X$ be a BIBD and let $f \in W(X)$. Then $f$ is an open map and in particular, the image of $f$ is a subsystem of $X$. Furthermore, there is an integer $d$ such that for every $y \in Im(f)$, $|f^{-1}(y)| =d$.

\et

Thus the partition induced by $f$ on $Dom(f)$, is a uniform partition: all classes have the same number of elements. We call the integer $d$ the degree of $f$ and write $d=deg(f)$. It follows that $|Dom(f)|=deg(f)|Im(f)|$. See \cite{bibdtrans} where a similar result for $BIBDs$ with arbitrary $\lambda$ is called the Homogeneous Lemma.

We have called a $PBD$ $X$ subsystem-free if the only subsystems of $X$ are the empty set, the points, the blocks and the whole set of points. That is, $X$ is subsystem-free if its only subsystems are the subsystems that every $PBD$ has. Equivalently, $X$ is subsystem-free if the corresponding matroid of $X$ has no proper extension to a larger $BRSC$ on the same point set.

The next theorem shows that subsystem-free $STSs$ on more than 9 points have small Wilson monoids. Thus the only subsystem-free $STS$ with non-small Wilson monoid are the Fano plane and the affine geometry $AG(2,3)$ as described above.

\bp \label{smallsts}

Let X be a subsystem-free STS on $v > 9$ points. Then $W(X)$ is a small monoid.

\ep

\proof
Let $f\in W(X)$. We have noted that $f$ is an open map and in particular, its range is a subsystem. Since $X$ is subsystem-free, the only possible ranges have size 0,1,3,$v$, where $v=|X|$. To prove that $W(X)$ is small, we must negate the possibility that the range of $f$ has 3 points, that is that the range of $f$ is a block.

 So assume that the range of $f$ is a block $b$ of $X$. Let $d=deg(f)$ as per Theorem \ref{uniform}. We recalled that $v$ is congruent to either 1 or 3 modulo 6 and we break up the proof into 2 cases.

1) $\bf v\equiv 1 \text{\bf  mod}(6)$

We know that $|Dom(f)|=d|Im(f)|=3d$ so that $Dom(f)$ is an open set with cardinality divisible by 3. The co-domain $D$, that is, the points on which $f$ is not defined is a subsystem. Let $y \in b$. Then we also have that
$f^{-w}(y)=D \cup f^{-1}(y)$ is a subsystem as well.

The open sets of $X$ have size 0,$v-3,v-1,v$. As noted above, $|Dom(f)|$ is a positive integer divisible by 3. Therefore, in this case, $|Dom(f)|=v-1$ and thus $d=\frac{v-1}{3}$. It follows that the subsystem $D \cup f^{-1}(y)$ has cardinality $\frac{v+2}{3}$. Given the possible sizes of subsystems, it follows that $v$ is at most 7, contradicting the assumption that $v>9$.

2) $\bf v\equiv 3 \text{\bf  mod}(6)$

Arguing as in the first case, $Dom(f)$ is an open set with cardinality divisible by 3. The possibilities are then either $|Dom(f)|=v-3$ or $|Dom(f)|=v$.

If $|Dom(f)|=v-3$, then $d=\frac{v-3}{3}$ and the subsystem $D \cup f^{-1}(y)$ has cardinality $\frac{v+6}{3}$. Given the possible sizes of subsystems, it follows that $v$ is at most 3, contradicting the assumption that $v>9$.

If $|Dom(f)|=v$, then $d=\frac{v}{3}$ and thus $D$ is empty and the subsystem $D \cup f^{-1}(y)$ has cardinality $\frac{v}{3}$. Again, it follows that $v$ is at most 9 and this is a contradiction. \qed

We now return to our survey of the Wilson monoids of $STS$.

{\bf v=13}

It is known that there are precisely 2 $STS$ up to isomorphism with $v=13$ \cite{tripsys}. It follows easily from Proposition \ref{subbound} and the congruential conditions on the orders of $STS$ that both of the $STS$ of order 13 are subsystem-free. Therefore each of them has a small Wilson monoid by Proposition \ref{smallsts}.

{\bf v=15}

Up to isomorphism there are 80 $STS$ of order 15 \cite[Pages 31-32]{handdesign}. Of these 57 are subsystem-free \cite[Table 1.29]{handdesign} and thus have small Wilson monoids by Proposition \ref{smallsts}. Among the 23 non-subsystem-free $STS$ of order 15 is the projective space of dimension 3 over the field of order 2, $P_{3,2}$. We have described its Wilson monoid in Example \ref{ProjWil}. In particular, it contains all $4 \times 4$ matrices over the field of order 2 as a submonoid and thus is not a small monoid.  The interested reader is welcome to survey the remaining 22 $STS$ of order 15.

{\bf v=19}

While the number of isomorphic $STS$ of order 15 was computed by hand in 1919 \cite{handdesign} it wasn't until the early 2000's that computer methods determined that there are 11,084,874,829 pairwise non-isomorphic $STS$ of order 19 \cite{sts19}. Of these, 10,997,902,498 are subsystem-free \cite{sts1921}. By Proposition \ref{smallsts} all have small Wilson monoids. Thus at least 99.2\% of the $STS$ of order 19 have small Wilson monoids.

We use the method of \cite{Wilson1, Wilson4} to construct an  $STS$ on 19 points with a non-small Wilson monoid. Let $L$ be a Latin Square on six points. We build a $GDD, \mathcal{G}(L)$ on the 18 points $R=\{R_{i}\mid i=1,...,6\} \cup C=\{C_{i}\mid i=1,...,6\} \cup X=\{1,2,3,4,5,6\}$. These three sets form the groups of $\mathcal{G}(L)$. The blocks are the triples $\{R_{i},C_{j},L(i,j)\}, i,j=1,...,6$. It is easy to see this forms a $GDD$ and more precisely a transversal design $TD(3,6)$. That is, the
$GDD$ has 3 groups each with 6 points. We have the associated $\{3,6\}-PBD$ by considering the groups to be blocks of order 6. We add a new point $p$ to $\mathcal{G}(L)$ and replace each of $R \cup \{p\}, L\cup \{p\}$ and $X \cup \{p\}$ by copies of the Fano plane. We now have defined an $STS, S(L)$ on 19 points that has 3 copies of the Fano plane that intersect pairwise in the point $p$. Let $Y=(\{x,y,z\},\{\{x,y,z\}\})$ be a trivial $PBD$ on 3 points. The function $f:S(L)\rightarrow Y$ with co-domain $\{p\}$ and that sends $R$ to x, $C$ to $y$ and $X$ to $z$ is a continuous map. Since any block of the form $\{R_{i},C_{j},L(i,j)\}$ is a section of $f$, it follows that $S(L)$ is of Wilson type and thus as mentioned previously, $W(S(L))$ is not a small monoid.

Of the 11,084,874,829 $STS$ of order 19, only 10,489, less than one in a million, are of Wilson type with a split idempotent of rank 3 and fibre a $TD(3,6)$ \cite{sts1921}. This paper also shows that there are precisely 2,156,186 $STS$ of Wilson type with a split idempotent of rank 3 and fibre a $TD(3,7)$. At the current time, there is no classification of all $STS$ of order 21. Current algorithms do not allow for a count of all $STS$ of order 21 in less than years of computer time. The
paper \cite{stssubs} determines the number of isomorphism classes of $STS$ on 21 points with a subsystem of order 9 and also those on 27 points with a subsystem of order 13.

An $STS$ is {\em rigid} if its automorphism group is trivial. Babai \cite{stsrigid} proved that almost all $STS$ are rigid. That is, the proportion of
such objects of an admissible order $n$ admitting non-trivial automorphisms tends to zero as $n \rightarrow \infty$. For n=19, of the 11,084,874,829 $STS$ up to isomorphism, only 164,078 have a non-trivial automorphism group.

Combining with the results on the number of $STS$ on 19 points \cite{sts19, sts1921}, there are 10,998,096,084 subsystem-free rigid $STS$ of order 19. All of their Wilson monoids are small monoids with trivial group of units. Thus more than 99\% of the $STS$ of order 19 have monoids consisting of a 0-simple semigroup with an identity element adjoined as Wilson monoids and thus the translational hull of the 0-simple semigroup is obtained by just adding an identity element.

As far as we know, there are no asymptotic results on subsystem-free $STS$. It seems reasonable for the results on $STS$ of order 19 that almost all $STS$ are subsystem-free. This would in turn mean that almost all Wilson monoids of $STS$ are small monoids with trivial group of units.

\section{Wilson Monoids of Pairwise Balanced Designs With One Block of Size Greater than 2}

In the previous section we looked in detail at the structure of Wilson monoids of Steiner Triple Systems. These are the smallest collection of $PBDs$, each of whose blocks has size greater than 2. In this section we look at the collection of $PBDs$ that have exactly one block of size greater than 2. For these, we can give the detailed structure of their Wilson monoids from the local (Green's relations) and global (various complexity functions) points of view.

Let $l \geq 3$ and $d \geq 1$ be given. Define $M(l,d)$ to be the $\{2,l\}-PBD$ with points $V=\{1, \ldots , l+d\}$. Let 
$L= \{1,...,l\}$. The blocks of $M(l,d)$ are $L$ together with $\{ \{ i,j\} \mid l+1 \leq j \leq l+d, 1 \leq i < j \}$.
That is, $M(l,d)$ has exactly one block of size greater than 2 and all the blocks of size 2 needed to ensure that we have a $PBD$. Notice that $M(l,1)$ is the Near Pencil of Example \ref{NP}. Let $W(l,d)$ be the Wilson monoid of $M(l,d)$. Let $D =\{l+1, \ldots , l+d\}$.

We begin by computing the subsystems and the open subsets of $M(l,d)$.

\bl \label{clopen}

Let $l \geq 3$, $d \geq 1$.

\begin{itemize}
  \item[(i)] {A subset $X$ of the points of $M(l,d)$ is a subsystem if and only if $L \subseteq X$ or $|L \cap X| \leq 1$.}

  \item[(ii)] {A subset $X$ of the points of $M(l,d)$ is open if and only if $X \subseteq D$ or $|L \cap X| \geq l-1$}.
\end{itemize}

\el

\proof
(i) It is easy to check that each such set $X$ is a subsystem of $M(l,d)$. Conversely, if $X$ is a subsystem that contains at least 2 points from $L$, then it contains $L$ by definition of a subsystem and the definition of $M(l,d)$.

(ii) follows from (i) by taking complements of sets. \qed

We now characterize the partial functions $f:V \rightarrow V$ that are in $W(l,d)$.

\bl \label{inW}

 Let $l \geq 3$, $d \geq 1$ and let  $f:V \rightarrow V$ be a partial function.

\begin{itemize}
  \item[(ii)] {If $L$ is not contained in $Dom(f)$, then $f \in W(l,d)$ if and only if $Dom(f)$ is open and $|f(L)| \leq 1$.}

  \item[(ii)] {If $L$ is contained in $Dom(f)$, then $f \in W(l,d)$ if and only if $|f(L)| = 1$ or $f$ restricted to $L$ is a permutation from $L$ to itself.}
\end{itemize}

\el

\proof
(i) Assume that $f \in W(l,d)$. Then it follows from Proposition \ref{morph1} that the domain of $f$ is open. Furthermore, since $L$ is not contained in the domain of $f$ it also follows from  Proposition \ref{morph1} that $|f(L)| \leq 1$. Conversely, let $f:V \rightarrow V$ be a partial function whose domain is open and does not contain $L$ and is such that $|f(L)| \leq 1$. Since the image under $f$ of any block of size 2 of $M(l,d)$ is either of size at most 1, another block of $M(l,d)$ or a 2 element subset of $L$ it follows from Proposition \ref{morph1} that $f \in W(l,d)$.

(ii) Assume that $L$ is contained in $Dom(f)$ and that $f \in W(l,d)$. It follows
from Proposition \ref{morph1} that if $|f(L)| > 1$, then $f$ restricted to $L$ is a permutation from $L$ to itself since $L$ is the unique block of $M(l,d)$ of size greater than 2. Conversely assume that $f$ is such that $|f(L)| = 1$ or $f$ restricted to $L$ is a permutation from $L$ to itself and that $L$ is contained in the domain of $f$. Then $Dom(f)$ is open by Lemma \ref{clopen}.  Furthermore, as in part (i), the image under $f$ of any block of size 2 of $M(l,d)$ is either of size at most 1, another block of $M(l,d)$ or a 2 element subset of $L$. It follows from Proposition \ref{morph1} that $f \in W(l,d)$. \qed

This Lemma allows us to characterize the possible images of elements of $W(l,d)$.

\bc
Let $l \geq 3$, $d \geq 1$ and let $X \subseteq V$.
Then there is an element $f \in W(l,d)$ with $Im(f)=X$ if and only if $L \subseteq X$ or $|X| \leq d+1$.

\ec

\proof
Assume that $f \in W(l,d)$. By Lemma \ref{inW} it follows that if $L$ is not a subset of $X$, then $|f(L)| \leq 1$ and thus $|Im(f)| \leq d+1$.

Conversely, if $L \subseteq X$, then the identity function restricted to $X$, $1|_X$ is a member of $W(l,d)$ by Lemma \ref{inW} and has range $X$. Let then $X$ be a
subset of $V$ with $|X| \leq d+1$ and let $Y=L \cap X$ and $Z = D \cap X$. If $Y$ is empty, then $X$ is the domain of $1_X$ and is open by Lemma \ref{clopen}. Thus $1_X \in W(l,d)$ by Lemma \ref{inW} and $X$ is the range of an element of $W(l,d)$.

Assume then that $|Y| = k$ and that $k>0$. If $|Z|=r$, then from $r+k \leq d+1$ it follows that $|D \setminus Z| \geq k-1$. Pick an element $j \in Y$ and a subset $W$ of $D \setminus Z$ with $|W| = k-1$ and let $g$ be a bijection from $W$ to $Y - \{j\}$. Then the partial function $f:V \rightarrow V$ with domain $L \cup Z \cup W$ defined by $f(v)=j$ if $v \in L$, $f(v)=g(v)$ if $v \in W$ and $f(v)=v$ if $v \in Z$ is in $W(l,d)$ by Lemma \ref{inW} and has image $X$. \qed

Recall that if $f:X\rightarrow X$ is a partial function, then the kernel of $f$ is the equivalence relation $Ker(f)$ on $Dom(f)$ defined by $(x,y) \in Ker(f)$ if and only
if $f(x)=f(y)$, $x,y \in Dom(f)$. We now characterize the kernels of elements of $W(l,d)$. If $\sim$ is an equivalence relation on a set $X$ and $Y \subseteq X$, then
the restriction $\sim|_Y$ is the equivalence relation on $Y$ defined by $\sim \cap (Y \times Y)$. In terms of partitions restriction to $Y$ has classes obtained by taking the non-empty intersections of classes of $\sim$ with $Y$.

\bl \label{Kernels}

Let $l \geq 3$, $d \geq 1$. Let $Y \subseteq V$ be an open set and let $\sim$ be an equivalence relation on $Y$. 

\begin{itemize}

  \item[(i)] {If $L \subseteq Y$, then there is an $f \in W(l,d)$ such that $Ker(f)=\; \sim$ if and only if either $\sim|_{L}$ is the identity relation on $L$ or $\sim|_{L}$ is the universal relation on $L$.}

  \item[(ii)] {If $L$ is not contained in $Y$ then there is an $f \in W(l,d)$ such that $Ker(f)=\; \sim$ if and only if  $\sim|_{L \cap Y}$ is the universal relation on
$L \cap Y$.}

\item[(iii)] {Let $f \in W(l,d)$. Then there is an idempotent $e \in W(l,d)$ such that $Ker(f)=Ker(e)$.}

\end{itemize}

\el

\proof
(i) Assume that $L \subseteq Y$. If $f \in W(l,d)$, then $\sim|_L$ is either the identity relation or the universal relation on $L$ by Lemma \ref{inW}.
Conversely, let $\sim$ be an equivalence relation on $Y$ such that $\sim|_L$ is the identity relation. We define a partial function $f$ with domain $Y$ by first having it be the identity function on $L$. Let $x \in Y \setminus L$. Then the $\sim$ equivalence class of $x$ contains at most one element $l(x)$ of $L$. For such classes, we extend the definition of $f$ so that $f(x)=l(x)$. The remaining equivalence classes of $\sim$ are contained in $Y \setminus L$. Let $Z$ be such a class.  Pick an element $z \in Z$ and extend the definition of $f$ by sending each element of $Z$ to $z$. Doing this for each such class defines a partial function $f$ with domain $Y$ and kernel $\sim$. Now $f \in W(l,d)$ by Lemma \ref{inW}. 


Assume that $\sim|_{L}= L \times L$. Let $Z$ be an equivalence class of $\sim$. Either $L \subseteq Z$ or $L \cap Z$ is the empty set. Pick a fixed element of each equivalence class. The partial function with domain $Y$ that sends an element to the representative of its class
is in $W(l,d)$ by Lemma \ref{inW} and has kernel equal to $\sim$. 

(ii) Now assume that $f \in W(l,d)$ and that $L$ is not contained in $Y$. Then by Lemma \ref{inW},  $\sim|_{L \cap Y}$ is the universal relation on
$L \cap Y$. Conversely, assume that $Y$ is an open set and that $\sim|_{L \cap Y}$ is the universal relation on
$L \cap Y$. Pick an element in each class of $\sim$. The function $f$ that sends each element of $Y$ to its representative is then an element of $W(l,d)$ by 
Lemma \ref{inW} and has kernel $\sim$. 

(iii) Note that the functions constructed in the proofs of parts (i) and (ii) are idempotents of $W(l,d)$. \qed


A semigroup $S$ is said to be a regular semigroup, if each of its elements is regular. Important regular semigroups include groups, inverse semigroups (defined by the property that each element has a unique inverse), the monoid of all functions (either total or partial) on a set and the monoid of all $n \times n$ matrices over a field. Despite these important examples, we now note that $W(l,d)$ is never a regular monoid.

\be

Let $l \geq 3$, $d \geq 1$. Let $i\neq j \in L$. Consider the total function $f:V \rightarrow V$ defined by $f(v)=i \text{ if } v \in L \text{ and } f(v)=j \text{ if } v \in D$.
Then $f \in W(l,d)$ by Lemma \ref{inW}. But $f$ is not a regular element of $W(l,d)$ by Corollary \ref{regimage}, since $Im(f)$ is not a subsystem.

\ee

\bq

We do not know of the existence of a Steiner triple system, $X$ such that its Wilson monoid $W(X)$ is not regular. As the preceding section showed, generically, Wilson monoids of Steiner triple systems seem to be small monoids, all of which are regular.

\eq

We now describe the regular elements of $W(l,d)$. We first recall some basic properties of Green's relations. See \cite{CP, qtheory} for more details. Let $M$ be a finite monoid.

Green's relations $\mathcal{R}$, $\mathcal{L}$ and $\mathcal{J}$ are defined on $M$ by
\begin{itemize}
\item $m\mathrel{\mathcal{L}} n$ if $Mm=Mn$;
\item $m\mathrel{\mathcal{R}} n$ if $mM=nM$;
\item $m\mathrel{\mathcal{J}} n$ if $MmM=MnM$.
\end{itemize}
The $\mathcal{L}$-class of $m\in M$ is denoted by $L_m$ and similar notation is used for $\mathcal{R}$- and $\mathcal{J}$-classes.  One defines the $\mathcal{L}$-order on $M$ by $m\leq_{\mathcal{L}} n$ if $Mm\subseteq Mn$.  The quasi-orders $\leq_{\mathcal{R}}$ and $\leq_{\mathcal{J}}$ are defined analogously.


The set of idempotents of $M$ is denoted by $E(M)$. Regularity of an element $m \in M$ is equivalent to each of the following: $L_m\cap E(M)\neq \emptyset$; $R_m\cap E(M)\neq \emptyset$; and $J_m\cap E(M)\neq \emptyset$ (the last equivalence uses finiteness).  A $\mathcal{J}$-class is called \emph{regular} if it contains an idempotent or, equivalently, contains only regular elements.  An important fact about finite monoids is that they enjoy a property called \emph{stability} which states that

$${xy\mathrel{\mathcal{J} x\iff xy\mathrel{\mathcal{R}} x \text{ and }
xy\mathrel{\mathcal{J} y\iff xy\mathrel{\mathcal{L}} y}}}$$

for $x,y\in M$~\cite[Theorem~1.13]{BenBook}.  One consequence of stability is that the intersection of any $\mathcal{R}$-class and $\mathcal{L}$-class in a $\mathcal{J}$-class is non-empty.  Another fact about finite semigroups that we shall use is that if $J$ is a $\mathcal{J}$-class such that $J^2\cap J\neq \emptyset$, then $J$ is regular (cf.~\cite[Corollary~1.24]{BenBook}).

Since Wilson monoids are explicitly given as submonoids of the monoid $PF(V)$ of all partial functions, we first quickly recall how to describe idempotents and Green's relations on this monoid. These results are classical and easy to prove. We have defined the kernel of a partial function as an equivalence relation on its domain, but we identify it with the corresponding partition on the domain. For a partial function $f:X \rightarrow X$, let $Fix(f)=\{x \in Dom(f)\mid f(x)=x\}$.

\bp \label{partialmon}

Let $V$ be a set and $f,g \in PF(V)$, the monoid of all partial functions on $V$.

\begin{itemize}
  \item[(i)] {$f$ is an idempotent if and only if $Im(f) = Fix(f)$.}
  \item[(ii)] {$f \mathcal{R} g$ if and only if $Im(f)= Im(g)$.}
  \item[(iii)] {$f \mathcal{L} g$ if and only if $Ker(f)= Ker(g)$.}
  \item[(iv)] {$f \mathcal{J} g$ if and only if $|Im(f)|=|Im(g)|$.}
\end{itemize}

\ep

The following result is also well known and we include it for completeness sake.

\bp \label{Greenids}

Let $M$ be a monoid and $N$ be a submonoid of $M$.

\begin{itemize}
  \item[(i)] {Let $e,f$ be idempotents in $M$. Then $e \mathcal{R} f \text{ if and only if } ef=f  \text{ and }fe=e$.}
  \item[(ii)] {Let $e,f$ be idempotents in $M$. Then $e \mathcal{L} f \text{ if and only if } ef=e  \text{ and }fe=f$.}
  \item[(iii)] {Let $x,y$ be regular elements of $N$. Then $x \mathcal{R} y$ in $N$ if and only if $x \mathcal{R} y$ in $M$. The dual statement for $\mathcal{L}$ also holds.}

\end{itemize}

\ep

\proof
Clearly if $ef=f$ and $fe=e$ then $e\mathcal{R}f$. Conversely, if there are elements $x,y$ in $M$ such that $ex=f \text{ and } fy = e$, then $ef=eex=ex=f$ and similarly $fe=e$. A dual proof works for $\mathcal{L}$. This proves (i) and (ii).

Now assume that $x,y$ are regular elements of $N$ and $x \mathcal{R} y$ in $M$. Since $x$ and $y$ are regular elements of $N$, there are idempotents $e,f$ in $N$, such that $x\mathcal{R}e$ and $y\mathcal{R}f$ in $N$. It follows that $e\mathcal{R}f$ in $M$. Therefore, $ef=f \text{ and }fe=e$ by part 1. These equations also hold in $N$ and thus $e\mathcal{R}f$ in $N$. We then have that $x\mathcal{R}e\mathcal{R}f\mathcal{R}y$ in $N$. A dual proof holds for $\mathcal{L}$. \qed

We now characterize which idempotents of $PF(V)$ belong to $W(l,d)$. The following follows immediately from Lemma \ref{inW}.

\bl \label{Wilsids}

Let $l \geq 3$, $d \geq 1$ and let $e$ be an idempotent in $PF(V)$. 

\begin{itemize}
  \item[(i)] {If $L \subseteq Dom(e)$, then $e \in W(l,d)$ if and only if $e|_L$ is the identity function on $L$} or $|e(L)|=1$.
	
  \item[(ii)] {If $L \nsubseteq Dom(e)$, then $e \in W(l,d)$ if and only if $Dom(e)$ is an open subset and $|e(L)| \leq 1$.}
\end{itemize}

\el




We can now describe the regular elements of $W(l,d)$.

\bl \label{regelts}

Let $l \geq 3$, $d \geq 1$. An element $f \in W(l,d)$  is regular if and only if $f|_L$ is a permutation of $L$ or $|Im(f) \cap L| \leq 1$. 


\el

\proof
First assume that $f \in W(l,d)$ and that $L \subseteq Im(f)$. In particular, $|Im(f) \cap L| > 1$. By Lemma \ref{inW}, either $f|_{L}$ is a permutation of $L$ or $|f(L)| \leq 1$.

If $f|_{L}$ is a permutation of $L$, for each element $x \in Im(f) - L$ pick an element $\bar{x}$ such that $f(\bar{x})=x$. Note that $\bar{x} \in D$. Define a partial function $g:V \rightarrow V$ with $Dom(g)= Im(f)$ by $g|_{L}=(f|_{L})^{-1}$ and for each $x \in Im(f) - L, g(x)=\bar{x}$. Then $g \in W(l,d)$ by Lemma \ref{inW}. It is routine to calculate that $fgf=f$ and thus $f$ is a regular element of $W(l,d)$.

Now let $f$ be such that $L \subseteq Im(f)$ and that $|f(L)| \leq 1$. Since $L \subseteq Im(f)$, there is a subset $X \subseteq D$ such that $|X|=l-1$ and $f(X) \subseteq L$ and also satisfies $|f(X)|=l-1$.
Let $g$ be any element of $W(l,d)$. It follows from Lemma \ref{inW} that either $g(f(X)) \subset L$ or that $|g(f(X))| \leq 1$. In both cases, it follows from our assumption that $|f(L)| \leq 1$ that $|fgf(X)| \leq 1$. Since $l-1 > 1$, it follows that $fgf \neq f$. This completes the proof in the case that $L \subseteq Im(f)$.


Now assume that $f \in W(l,d)$ is a regular element and that $L$ is not a subset of $Im(f)$. Since $f$ is a regular element, it follows from Corollary \ref{regimage} that $Im(f)$ is a subsystem of $M(l,d)$. It follows from Lemma \ref{clopen} that $|Im(f) \cap L| \leq 1$.

Conversely assume that $f \in W(l,d)$ has $|Im(f) \cap L| \leq 1$. It follows from Lemma \ref{inW} that $|f(L)| \leq 1$ as well. We have two cases.

{\bf Case 1}: $Im(f) \cap L = \emptyset$. 

$Im(f) \subseteq D$ and is thus an open set by Lemma \ref{clopen}. For each $x \in Im(f)$ pick an $\bar{x} \in Dom(f)$ such that $f(\bar{x})=x$. Define $g:V \rightarrow V$ to be the function such that $Dom(g)=Im(f)$ and with $g(x)=\bar{x}$. It is clear that $fgf=f$ and $g \in W(l,d)$ by Lemma \ref{inW} and $fgf=f$.  





{\bf Case 2} $Im(f) \cap L = \{j\}$ for some $j \in L$.

For each $x \in Im(f)$ pick $\bar{x}$ such that $f(\bar{x}) = x$. Define $g:L \cup Im(f) \to W(l,d)$ by $g(x) = \bar{x}$ if $x \in Im(f)$ and $g(x) = \bar{j}$ if $x \in L$. Then $g \in W(l,d)$ by Lemma \ref{inW} and $fgf=f$. \qed

\brm

It follows from 
Corollary \ref{regimage} 
that the image of any regular element of $W(l,d)$ is a subsystem of $M(l,d)$. Moreover if $X$ is a subsystem, then there is a regular element $f$ such that $Im(f)=X$. Indeed, if $X$ is a subset of $D$ or $L \subseteq X$, then the identity function restricted to $X$ is in $W(l,d)$ by Lemma \ref{inW} and has image $X$.  If $L \cap X =\{j\}$ for some $j \in L$, then the partial function with domain $L \cup X$ that sends all of $L$ to $j$ and is the identity of 
$X \setminus \{j\}$ is an idempotent in $W(l,d)$ by Lemma \ref{inW} and has range $X$. 

Despite this, it does not follow that every element of $W(l,d)$ with image a subsystem is a regular element. For example, if $l=d=3$, the partial function $f$ with domain $\{4,5,6\}$ and such that $f(4)=1,f(5)=2,f(6)=3$ is in $W(l,d)$ by Lemma \ref{inW}, has range the subsystem $L$, but is not regular by Lemma \ref{regelts}.

\erm

We now describe the $\mathcal{J}$ relation for regular elements in $W(l,d)$. We divide this into two cases, depending on whether $L$ is or is not a subset of the
image of a regular element. We first look at the case when $L\nsubseteq Im(f)$, so that $|Im(f) \cap L| \leq 1$ by Lemma \ref{regelts}. In this case, it follows from Lemma \ref{inW} that $|f(L)| \leq 1$ as well. Let $0 \leq i\leq d+1$ and define 
$$J_{i}=\{f \in W(l,d) \; \big{\vert} \; |Im(f) \cap L| \leq 1 \text{ and }|Im(f)|=i\}.$$

Recall that an ideal $I$ of a monoid $M$ is said to be prime if its complement $M - I$ is
a submonoid or, equivalently, $I$ is a proper ideal and $ab \in I \text{ implies that } a \in I \text{ or } b \in I$.

\bl \label{constantJ}

Let $l \geq 3$, $d \geq 1$.

\begin{itemize}
  \item[(i)] {Let $I=\{f \in W(l,d) \; \big{\vert} \; |f(L)| \leq 1\}$. Then $I$ is a prime ideal of $W(l,d)$.}

  \item[(ii)] {Let $f,g \in I$ be regular elements of $W(l,d)$. Then $f\mathcal{J}g$ if and only if $|Im(f)|=|Im(g)|$.}
  \item[(iii)] {The regular $\mathcal{J}$-classes contained in $I$ are precisely $\{J_{i} \mid 0 \leq i\leq d+1\}$.}
  \item[(iv)] {The unique maximal $\mathcal{J}$-class of $I$ is $J_{d+1}$.}
\end{itemize}

\el

\proof
(i) Let $f \in I$ and let $k,h \in W(l,d)$. By Lemma \ref{inW}, either $h(L)=L$ or $|h(L)| \leq 1$. In both cases, $f \in I$ implies that $|fh(L)|\leq 1$ and thus
$|kfh(L)| \leq 1$. Therefore, $I$ is an ideal of $W(l,d)$. Since every element of $W(l,d)$ either restricts to a permutation on $L$ or is such that its image has rank at most 1 on $L$ by Lemma \ref{inW}, it immediately follows that $I$ is a prime ideal of $M$.

(ii) Let $f,g$ be regular elements in $I$. If $f\mathcal{J}g$, then $|Im(f)|=|Im(g)|$ by Proposition \ref{partialmon}. Conversely, assume that $|Im(f)|=|Im(g)|$ for
$f,g$ regular elements in $I$.

We first consider the case that $|Im(f)|=d+1$. It follows from Lemma \ref{regelts} that $Im(f) = D \cup \{j\}$ and $Im(g) = D \cup \{j'\}$ for some $j,j' \in L$. Let $h:V\rightarrow V$ be a permutation that maps $L$ onto $L$ and such that $h(j)=j'$. Then $h$ belongs to the group of units of $W(l,d)$ and thus $hf\mathcal{L}f$. Since $Im(hf)=Im(g)$, $hf\mathcal{R}g$ by Proposition \ref{partialmon} and Proposition \ref{Greenids}(iii). Therefore, $f\mathcal{J}g$.

Assume then that $|Im(f)|=|Im(g)| \leq d$. Assume that there is a $j \in L$ that belongs to $Im(f)$. Let $v \in D$ be an element not in the image of $f$ and let $h:V\rightarrow V$ be the identity on $D - \{v\}$, send $L$ to $v$ and $v$ to $j$.
Then $h \in W(l,d)$ by Lemma \ref{inW}. Then $Im(hf) \subseteq D$ and since $h^{2}f=f$, $hf \mathcal{L} f$. Therefore, we can assume that both $Im(f) \text{ and }Im(g)$ are contained in $D$.

Let $k:D\rightarrow D$ be a permutation such that $k(Im(f))=Im(g)$. Then considered as partial functions on $V$ both $k,k^{-1}$ belong to $W(l,d)$ and since $k^{-1}kf=f$, we have that $kf\mathcal{L}f$. Since $Im(kf)=Im(g)$, we have $kf\mathcal{R}g$ by Proposition \ref{partialmon} and Proposition \ref{Greenids}(iii). Therefore, $f\mathcal{J}g$.

(iii) This is an immediate consequence of (ii).

(iv) Let $f \in I$. Then $|Im(f)|\leq d+1$. If $f(L)= \emptyset$, then $Dom(f) \subseteq D$. Fix an element $j \in L$. The total function $g:V\rightarrow V$ such that $g(L)=\{j\}$ and is the identity on $D$ is in $J_{d+1}$ and $fg=f$. Therefore $f \leq_{\mathcal{J}}g$. If $f(L)=\{v\}, v \in V$, we choose  $j\in Dom(f) \cap L$ and define a total function $g$ $g$ and such that $g(x)=j$ if $x \in L$ and $g(x)=x$ if $x \in D$. Again, $fg=f$ and $g \in J_{d+1}$. Therefore $f \leq_{\mathcal{J}}g$ in this case as well. \qed



We now turn to the description of $\mathcal{J}$-classes for regular elements $f \in W(l,d)$ such that $L \subseteq Im(f)$. Let $0 \leq i \leq d$. Define 
$$J_{L,i}=\{f\in W(l,d) \mid f \text{ is a regular element}, L \subseteq Im(f), |Im(f)|=l+i\}.$$
Notice that
$J_{L,d}$ is the group of units of $W(l,d)$.

\bl \label{permutationJ}

Let $l \geq 3$, $d \geq 1$.
\begin{itemize}
\item[(i)] Let $N = \{f \in W(l,d)\mid f|_{L} \text{ is a permutation on L}\}$. Then $N$ is a regular submonoid of $W(l,d)$ and is a union of $\mathcal{J}$-classes of $W(l,d)$.
\item[(ii)] Let $f \in W(l,d)$ be a regular element with $L \subseteq Im(f)$. Then $f \in N$ and the $\mathcal{J}$-class of $f$ is $J_{L,i}$ where $|Im(f)|=l+i$.
\end{itemize}
\el

\proof
(i) By Lemma \ref{inW}, $N=W(l,d)-I$ where $I$ is the ideal discussed in Lemma \ref{constantJ}. Since $I$ is a prime ideal, $N$ is a submonoid of $W(l,d)$. Furthermore, every element of $N$ is regular by Lemma \ref{regelts}. Finally it follows from \cite[Lemma 2.2]{BenBook} that $N$ is a union of $\mathcal{J}$-classes of $W(l,d)$.

(ii) Let $f$ be a regular element with $L \subseteq Im(f)$. It follows from Lemma \ref{regelts} that $f \in N$. Let $J$ be the $\mathcal{J}$-class of $f$ and let $g \in J$. Then $g$ is a regular element and $|Im(g)|  =l+i$, since $f\mathcal{J}g$ as elements of $PF(V)$. By part (i), 
$L \subseteq Im(g)$ and thus $g \in J_{L,i}$. Conversely, let $g \in J_{L,i}$. Then $Im(f)=L \cup X$ and $Im(g)=L \cup Y$ where $X$ and $Y$ are subsets of $D$ with $|X|=|Y|=i$. Let $h\in Sym(V)$ be any permutation that is the identity on $L$ and maps $X$ onto $Y$. Then $h$ is in the group of units of $W(l,d)$ and thus $hf\mathcal{L}f$. Since $Im(hf)=Im(g)$, it follows from Proposition \ref{partialmon} and Proposition \ref{Greenids}(iii) that $hf\mathcal{R}g$ and thus $f\mathcal{J}g$. \qed

We summarize the results on Green's relations for regular elements in $W(l,d)$. We use the notation from the previous lemmas.

\bt \label{GreenW}
Let $l \geq 3$, $d \geq 1$ and let $f$ and $g$ be regular elements of $W(l,d)$.
\begin{itemize}
\item[(i)]
$f\mathcal{R}g$ if and only if $Im(f)=Im(g)$.
\item[(ii)]
$f\mathcal{L}g$ if and only if $Ker(f)=Ker(g)$.
\item[(iii)]
There are precisely $2d+3$ regular $\mathcal{J}$-classes and they are $\{J_{i} \mid i=0,\ldots d+1\} \cup \{J_{L,i}\mid i=0,\ldots d\}$.
\item[(iv)]
The maximal subgroup of $J_{i}$ is the symmetric group $S_i$ on $i$ elements. The maximal subgroup of $J_{L,i}$ is $S_{l} \times S_i$.
\item[(v)]
Let $\Omega(l,d)$ be the poset of regular $\mathcal{J}$-classes of $W(l,d)$. Then $\{J_{i}\mid i=0,\ldots d+1\}$ and $\{J_{L,i}\mid i=0,\ldots d\}$ form chains of length $d+2$ and $d+1$ respectively in $\Omega(l,d)$.
\item[(vi)]
In $\Omega(l,d)$, $J_{L,i}$ covers precisely $J_{L,i-1}$ and $J_{i+1}$, for $i=0,\ldots d$.
\item[(vii)]
In $\Omega(l,d)$, $J_{i}$ covers precisely $J_{i-1}$ for $i=1,\ldots d+1$ and $J_0$ is the unique minimal element.
\item[(viii)]
The $\mathcal{J}$-classes $J_{L,d-1} \text{ and } J_{d+1}$ are the unique two maximal $\mathcal{J}$-classes less than the group of units in the poset of all $\mathcal{J}$-classes of $W(l,d)$.
\end{itemize}
\et

\proof
(i), (ii) and (iii) follow from Proposition \ref{partialmon}, Proposition \ref{Greenids}, Lemma \ref{constantJ} and Lemma \ref{permutationJ}.

It is well known that the maximal subgroup in a $\mathcal{J}$-class, $J$ in a monoid $M$ is isomorphic to the group of units of the monoid $eMe$ for any idempotent $e \in J$ \cite{qtheory}. Let $X$ be a subset of $D$ with $|X|=i$. Then $1_X$, the identity function restricted to $X$ is an idempotent that belongs to $J_i$. Every permutation of $X$ considered as a partial function on $V$ is in the group of units of $1_{X}W(l,d)1_X$ and thus the maximal subgroup of $J_i$ is $S_i$, $0 \leq i \leq d$. 

We consider now the case of $J_{d+1}$.  Let $j \in L$. The total function $e$ that sends all of $L$ to $j$ and is the identity on $D$ is an idempotent in $W(l,d)$ by Lemma \ref{inW} and belongs to $J_{d+1}$. Let $\sigma$ be a permutation of $D \cup \{j\}$. Extend $\sigma$ to a total function $\bar{\sigma}$ on $V$ by letting $\bar{\sigma}$ agree with $\sigma$ on $D$ and by sending each element of $L$ to $\sigma(j)$. It is easy to check that $\bar{\sigma}$ is in the group of units $G$ of $eW(l,d)e$ and that the assignment of $\sigma$ to $\bar{\sigma}$ is an isomorphism of $S_{d+1}$ onto $G$, since every element of $G$ restricts to a permutation of $D \cup \{j\}$. Therefore, $G$ is isomorphic to $S_{d+1}$.

Similarly, $1_{L \cup X}$ the identity restricted to $L \cup X$ belongs to $J_{L,i}$. By Lemma \ref{inW} the invertible elements of $1_{L \cup X}W(l,d)1_{L \cup X}$ are precisely the permutations of $V$ that restrict to permutations of both $L$ and $X$. Clearly, this group is isomorphic to $S_{l}\times S_{i}$. This proves (iv).

Fix a chain of subsets $X_{0}=\emptyset \subset X_{1} \subset X_{2} \subset \ldots \subset X_{d}=D$ of subsets of $D$ with $|X_{i}|=i$. Recall that the collection of idempotents of a monoid $M$ is partially ordered by declaring that for $e,f \in E(M), e \leq f \text{ if and only if } e=ef=fe$, equivalently, that $f$ is below $e$ in both the $\mathcal{R}$ and the $\mathcal{L}$ orders of $M$. Furthermore, for regular $\mathcal{J}$-classes, $J,J'$ of $M$, $J \leq_{\mathcal{J}} J'$ if and only if there are idempotents $e \in J, f \in J'$ with $e \leq f$. Moreover, if $J \leq J'$ then for each idempotent $f \in J'$, there is an idempotent $e \in J$ such that $e \leq f$ \cite{qtheory}.

For each set $X_{i}, 0 \leq i \leq d+1$ and for each set $L \cup X_{i}$ the identity function restricted to these sets is an idempotent in $W(l,d)$. They clearly form chains in the idempotent ordering and thus the $\mathcal{J}$ ordering, proving the assertion in (v). The proof of (vi) and (vii) follow from consideration of the idempotents defined here as well.

We turn to the proof of (viii). By Lemma \ref{constantJ}, $J_{d+1}$ is the unique maximal $\mathcal{J}$-class of $W(l,d)$ in the ideal $I$. Since $I$ is a prime ideal, any $\mathcal{J}$-class above $J_{d+1}$ must belong to the regular submonoid $N$ and thus be equal to $J_{L,i}$ for some $0 \leq i \leq d$ by Lemma \ref{permutationJ}. If $i <d$ and using the notation from the preceding paragraphs, the identity function restricted to $L \cup X_{i}$ is an idempotent $f$ that belongs to $J_{L,i}$. Every idempotent $e$ in $J_{d+1}$ has $D$ as a subset of its range. But $X_{i}$ is a proper subset of $D$ if $i < d$ and thus no idempotent in $J_{d+1}$ is below $f$ in the $\mathcal{R}$ order of $W(l,d)$ and thus there is no idempotent $e\in J_{d+1}$ such that $e \leq f$. Therefore $J_{d+1}$ is a maximal $\mathcal{J}$-class in the poset of all $\mathcal{J}$-classes of $W(l,d)$.

Finally, $J_{L,d-1}$ is covered by the group of units $J_{L,d}$ by part 6. of this Theorem and no $\mathcal{J}$ class in  $I$ can be above $J_{L,d-1}$ in the $\mathcal{J}$ order because $I$ is an ideal and $J_{L,d-1}$ belongs to the complement $W(l,d) - I$. This completes the proof of (viii) and of the theorem. \qed

The final topic of this section determines the complexity of $W(l,d)$ in the sense of Krohn-Rhodes decomposition theory \cite[Part II]{qtheory}. We recall some basic definitions and results. It follows from the Krohn-Rhodes Decomposition Theorem \cite[Theorem 4.1.30]{qtheory} that if $S$ is a finite semigroup, then $S$ divides, that is, $S$ is the homomorphic image of a subsemigroup of an iterated wreath product of finite groups and finite semigroups all of whose maximal subgroups are trivial. The least number of non-trivial groups in any such decomposition is called the (Krohn-Rhodes) complexity of $S$. We write $Sc$ for the complexity of $S$. The reverse complexity of $S$, denoted by $Sc^*$ is the complexity of the reverse semigroup $S^{op}$ of $S$. There are examples where the complexity and reverse complexity of a semigroup can differ by an arbitrary amount \cite[Chapters 7-9]{Arbib}. We summarize here the results from complexity theory that we need here. Some of these results are easy to prove and some require some of the deepest results of complexity theory.

\bt \label{comp}
\begin{itemize}
\item[(i)]
Let $S$ and $T$ be finite semigroups. If $S$ divides $T$, then $Sc \leq Tc$ and if $S$ and $T$ are non-empty, then $(S \times T)c=max\{Sc,Tc\}$.
\item[(ii)]
Let $V$ be a set. Then the complexity of the full transformation monoid and the monoid of all partial transformations on $V$ is $|V|-1$.
\item[(iii)]
Let $I$ be an ideal of a semigroup $S$. Then $Sc \leq (S/I)c + Ic$.
\item[(iv)]
Let $S$ be a semigroup that is equal to $SeS$ for some idempotent $e \in S$. Then $Sc = (eSe)c$.
\item[(v)]
Let $M$ be a small monoid. If the idempotent generated submonoid of $M$ has trivial subgroups then $Mc=Mc^{*}$ and is 0 if all subgroups of $M$ are trivial and equal to 1 otherwise.
\end{itemize}
\et

\proof
(i) These are well known facts about complexity \cite[Chapter 4]{qtheory}. 

(ii) This fact was first proved in \cite{resultsonfinite} for the full transformation semigroup. The results for the monoid of all partial functions can be proved similarly. 

(iii) This statement is the Ideal Theorem \cite{TilsonXII}, \cite[Theorem 4.9.17]{qtheory}. 

(iv) This statement follows from the Reduction Theorem \cite{TilsonXI},\cite[Theorem 4.9.16]{qtheory}. 

(v) This statement follows from Tilson's 2-$\mathcal{J}$-class Theorem \cite[Section 4.15]{qtheory}. \qed

We will now compute the complexity and reverse complexity of $W(l,d)$.  We need a few technical lemmas. We use the notation for $\mathcal{J}$-classes of $W(l,d)$ used above.

\bl \label{smallc}
Let $l \geq 3$, $d \geq 1$.
Let $K$ be the ideal of $W(l,d)$ generated by $J_{L,d-1}$. Then $W(l,d)/K$ is a small monoid with 0-minimal ideal the principal factor corresponding to $J_{d+1}$. Furthermore, $(W(l,d)/K)c = (W(l,d)/K)c^{*}=1$.

\el

\proof
By Theorem \ref{GreenW}(vi) and (viii), in order to prove that $W(l,d)/K$ is a small monoid with 0-minimal ideal the principal factor corresponding to $J_{d+1}$ it is enough to prove that every non-regular element of $W(l,d)$ is contained in $K$. By Lemma \ref{inW} and Lemma \ref{regelts} if $f$ is a non-regular element of $W(l,d)$, then $|Im(f)| \leq d+1$ and 
$|Im(f) \cap L| \geq 2$. Therefore, there is an $i \in D$ that is not in the image of $f$. The idempotent $e$ that is the identity function restricted to $V -\{i\}$ belongs to $J_{L,d-1}$. Since $ef=f, f \in K$ as desired.

We have proved that $W(l,d)/K=G \cup J_{d+1} \cup\{0\}$ is a small monoid, where $G$ is the group of units of $W(l,d)$. We claim that the idempotent generated submonoid of this small monoid has trivial subgroups and the second statement will then follow from Theorem \ref{comp}(v).

Let $e$ be an idempotent of $J_{d+1}$. Then $Im(e)=\{j\} \cup D$ for some $j \in L$. Since $e$ is an idempotent, the restriction of $e$ to $D$ is the identity function on $D$. Therefore if $f=e_{1}e_{2} \ldots e_{k}$ is the product of the idempotents $e_{i}, i=1,\ldots k$ in $J_{d+1}$, then $f$ restricted to $D$ is the identity on $D$. If $Im(f)=D$, then $f=0$ in $W(l,d)/K$. Otherwise, $|Im(f)|=\{j'\} \cup D$ for some $j' \in L$ and $f^{2}=f$ if $j'\in Dom(f)$ and $f^{2}=0$ otherwise, in $W(l,d)/K$. Therefore, in all cases, $f^{2}=f^{3}$ and all subgroups of the idempotent generated submonoid of $W(l,d)/K$ are trivial. This completes the proof. \qed

The next two lemmas will allow us to use induction in the proof of the main theorem.

\bl \label{K}
Let $l \geq 3$, $d > 1$.
Let $K$ be the ideal defined in the previous lemma. Let $e$ be the identity function restricted to $V-\{l+d\}$. Then $e\in J_{L,d-1}$ 
and $eKe$ is isomorphic to $W(l,d-1)$.

\el

\proof
It is clear from Lemma \ref{inW} that $e\in J_{L,d-1}$. Furthermore, if $f \in K$, then $efe =f$ if and only if both $Dom(f)$ and $Im(f)$ are contained in $V-\{l+d\}$. Therefore, by Lemma \ref{inW}, we can consider each such function to be a member of $W(l,d-1)$ and this assignment is easily seen to be an isomorphism between $eKe$ and 
$W(l,d-1)$.


\bl \label{one}
$W(l,1)c=W(l,1)c^{*}=1$.
\el

\proof
Let $S=W(l,1)$ and let $G$ be the group of units of $S$. By Theorem \ref{GreenW}, the complement of $G \cup J_{L,0} \cup J_{2}$ is an ideal $I$ of $S$. The only regular $\mathcal{J}$-classes in $I$ are $J_{1}$ and $J_{0}$. Therefore all subgroups of $I$ are trivial and therefore $Ic=Ic^{*}=0$. It follows from Theorem \ref{comp} that 
$Sc=(S/I)c \text{ and }Sc^{*}=(S/I)c^{*}$.

Let $T=S/I$. By Theorem \ref{GreenW} (viii),  both $J_{L,0} \cup \{0\}$ and $J_{2} \cup \{0\}$ are 0-minimal ideals of $T$ and $T$ is the union of these two ideals and its group of units $G$. Therefore, $T$ is a subdirect product of the Rees quotients $T / {(J_{L,0} \cup \{0\})}$ and $T / {(J_{2} \cup \{0\})}$. By Theorem \ref{comp} (i) it suffices to prove that each of these quotients has (reverse) complexity equal to 1.


We have seen in the last paragraph of the proof of Lemma \ref{smallc} that the idempotent generated subsemigroup of $J_{2} \cup \{0\}$ is aperiodic. Thus the (reverse) complexity of the small monoid $T / {(J_{L,0} \cup \{0\})}$ is equal to 1 by Theorem \ref{comp} (v).

Finally $T/(J_{2} \cup \{0\})$ is a small monoid with 0-minimal ideal $J_{L,0} \cup \{0\}$ and thus by Theorem \ref{comp}(v), it suffices to show that the idempotent generated subsemigroup of $J_{L,0} \cup \{0\}$ has trivial subgroups. If $f \in J_{L,0}$, then $|Im(f)|=|L|$ and since $L \subseteq Im(f)$ due to $f|_L$ being a permutation on $L$, we have that $Im(f)=L$ for all $f \in J_{L,0}$. Therefore, all elements of $J_{L,0}$ are in the same $\mathcal{R}$-class by Theorem \ref{GreenW}. It follows easily from Proposition \ref{Greenids} that the idempotents of $J_{L,0}$ form a right-zero semigroup and this proves the result. \qed

The next lemma gives a lower bound to the complexity functions of $W(l,d)$.

\bl \label{full}

Let $l \geq 3$, $d \geq 1$. The full transformation monoid on $d+1$ points is isomorphic to a subsemigroup of $W(l,d)$. Therefore $d \leq W(l,d)c$ and $d \leq W(l,d)c^{*}$.

\el

\proof
Let $j \in L$ and $f:\{j\} \cup D \rightarrow \{j\} \cup D$ be a total function. Then the total function $\bar{f}:V \rightarrow V$ defined by sending all elements of $L$ to $f(j)$ and if $d \in D$, then $\bar{f}(d)=f(d)$ belongs to $W(l,d)$. It is easy to see that the assignment $f$ to $\bar{f}$ is an injective morphism. The inequality for complexity and reverse complexity follow from Theorem \ref{comp}. \qed

Here is the main result on the complexity of $W(l,d)$.

\bt \label{Wilsoncomp}

Let $l \geq 3$, $d \geq 1$. Then the complexity and the reverse complexity of $W(l,d)$ are equal to $d$.

\et

\proof
We prove this by induction on $d$. The case of $d=1$ is given by Lemma \ref{one}. By Lemma \ref{full} we need only prove that $W(l,d)c \leq d$ and $W(l,d)c^{*} \leq d$.

Assume that the complexity and the reverse complexity of $W(l,d-1)$ are at most $d-1$, $d > 1$. Consider the ideal $K$ in $W(l,d)$ defined in Lemma \ref{smallc}. By Theorem \ref{comp}(iii), $W(l,d)c \leq ((W(l,d)/K)c+ Kc)$. Let $e$ be the identity function restricted to
$V-\{l+d\}$.  We claim that $K = KeK$.

Since $e \in J_{L,d-1} \subseteq K$, we have $KeK \subseteq K$. Conversely, let $f \in K$. If $l+d \notin Im(f)$, we immediately get $f = eef \in KeK$, hence we may assume that $l+d \in Im(f)$. Suppose first that $D \not\subseteq Im(f)$. Let $v \in D - Im(f)$ and let $g$ be the permutation of $V$ that exchanges $v$ and $l+d$ and is the identity elsewhere. Then $g \in W(l,d)$ and $f = (ge)e(gf) \in KeK$. Suppose now that $D \subseteq Im(f)$. In the proof of Lemma \ref{smallc}, it was noted that the non-regular elements of $W(l,d)$ satisfy $|Im(f)| \leq d+1$ and 
$|Im(f) \cap L| \geq 2$, hence $f$ must be regular and so $f \in J_d$ in view of Theorem \ref{GreenW}. Thus $Im(f) = D$.
Let $j \in L$ and let $g$ be the permutation of $V$ that sends $L$ to $l+d$, $l+d$ to $j$ and is the identity elsewhere. Then $g \in W(l,d)$ and $f = (ge)e(gf) \in KeK$ also in this case.

From Lemma \ref{K} and Theorem \ref{comp}(iv), the fact that $K =KeK$ and the inductive hypothesis, we get $Kc=W(l,d-1)c \leq d-1$. By Lemma \ref{smallc}, $(W(l,d)/K)c \leq 1$. This proves the result for complexity. A similar proof proves the result for reverse complexity as well. \qed

\section*{Acknowledgments}

The first author acknowledges support from the Binational Science Foundation (BSF) of the United States and Israel, grant number 2012080. The second author acknowledges support from the Simons Foundation (Simons Travel Grant Number 313548).
The third author was partially supported by CMUP (UID/MAT/00144/2019), which is funded by FCT (Portugal) with national (MCTES) and European structural funds through the programs FEDER, under the partnership agreement PT2020.

\bibliography{stubib}
\bibliographystyle{abbrv}

\vspace{1cm}

{\sc Stuart Margolis, Department of Mathematics, Bar Ilan University,
  52900 Ramat Gan, Israel}

{\em E-mail address:} margolis@math.biu.ac.il

\bigskip

{\sc John Rhodes, Department of Mathematics, University of California,
  Berkeley, California 94720, U.S.A.}

{\em E-mail addresses}: rhodes@math.berkeley.edu, BlvdBastille@gmail.com

\bigskip

{\sc Pedro V. Silva, Centro de
Matem\'{a}tica, Faculdade de Ci\^{e}ncias, Universidade do
Porto, R. Campo Alegre 687, 4169-007 Porto, Portugal}

{\em E-mail address}: pvsilva@fc.up.pt

\end{document}